\documentclass[12pt]{article}

\usepackage{epsfig}
\usepackage{amssymb, enumerate}

\newtheorem{theorem}{Theorem}
\newtheorem{definition}[theorem]{Definition}
\newtheorem{lemma}[theorem]{Lemma}
\newtheorem{proposition}[theorem]{Proposition}
\newtheorem{corollary}[theorem]{Corollary}
\newtheorem{example}[theorem]{Example}

\newtheorem{remark}[theorem]{\sc Remark}

\newcommand{\bbox}{\ \hfill\rule[-1mm]{2mm}{3.2mm}}

\title{{\sc Generalized Hopf Bifurcation for planar vector fields via the inverse integrating factor}\thanks{The
authors are partially supported by a DGICYT grant number
MTM2005-06098-C02-02.}}

\author{{\sc Isaac A. Garc\' \i a$^{\ (1)}$, H\'ector Giacomini$^{\ (2)}$ \& Maite Grau$^{\ (1)}$}}
\date{}

\begin{document}
\maketitle \vspace{-1.0cm}

\begin{abstract}
In this paper we study the maximum number of limit cycles that can
bifurcate from a focus singular point $p_0$ of an analytic,
autonomous differential system in the real plane under an analytic
perturbation. We consider $p_0$ being a focus singular point of
the following three types: non-degenerate, degenerate without
characteristic directions and nilpotent. In a neighborhood of
$p_0$ the differential system can always be brought, by means of a
change to (generalized) polar coordinates $(r, \theta)$, to an
equation over a cylinder in which the singular point $p_0$
corresponds to a limit cycle $\gamma_0$. This equation over the
cylinder always has an inverse integrating factor which is smooth
and non--flat in $r$ in a neighborhood of $\gamma_0$. We define
the notion of vanishing multiplicity of the inverse integrating
factor over $\gamma_0$. This vanishing multiplicity determines the
maximum number of limit cycles that bifurcate from the singular
point $p_0$ in the non-degenerate case and a lower bound for the
cyclicity otherwise. \par Moreover, we prove the existence of an
inverse integrating factor in a neighborhood of many types of
singular points, namely for the three types of focus considered in
the previous paragraph and for any isolated singular point with at
least one non-zero eigenvalue.
\end{abstract}
{\small{\noindent 2000 {\it AMS Subject Classification:} 37G15, 37G10, 34C07   }}

\section{Introduction and statement of the results \label{sect1} }

Let us consider a planar, real, analytic, autonomous differential
system with a singular point which we assume to be at the origin,
that is, we consider a differential system of the form:
\begin{equation}
\dot{x} \, = \, P(x,y), \qquad \dot{y} \, = \, Q(x,y), \label{eq1}
\end{equation}
where $P(x,y)$ and $Q(x,y)$ are real analytic functions in a
neighborhood $\mathcal{U}$ of the origin such that
$P(0,0)=Q(0,0)=0$. As usual, we associate to system (\ref{eq1})
the vector field $\mathcal{X}_0 \, = \, P(x,y) \partial_x \, + \,
Q(x,y) \partial_y$. We assume that the origin $p_0=(0,0)$ is an
isolated singular point, that is, there exists a neighborhood of
it without any other singular point, and we assume that it is a
monodromic singular point. Therefore, it is either a center (i.e.
it has a neighborhood filled with periodic orbits) or a focus
(i.e. it has a neighborhood where all the orbits spiral in forward
or in backward time to the origin).
\par We consider an analytic perturbation of system (\ref{eq1}) of
the form:
\begin{equation}
\dot{x} \, = \, P(x,y) \, + \, \bar{P}(x,y,\varepsilon), \qquad
\dot{y} \, = \, Q(x,y) \, + \, \bar{Q}(x,y,\varepsilon),
\label{eq2}
\end{equation} where $\varepsilon \in \mathbb{R}^p$ is the perturbation parameter,
$0<\|\varepsilon\|<<1$ and the functions
$\bar{P}(x,y,\varepsilon)$ and $\bar{Q}(x,y,\varepsilon)$ are
analytic for $(x,y) \in \mathcal{U}$, analytic in a neighborhood
of $\varepsilon=0$ and $\bar{P}(x,y,0)= \bar{Q}(x,y,0)\equiv 0$.
We associate to this perturbed system (\ref{eq2}) the vector field
$\mathcal{X}_\varepsilon \, = \, (P(x,y) \, + \,
\bar{P}(x,y,\varepsilon))\partial_x \, + \, (Q(x,y) \, + \,
\bar{Q}(x,y,\varepsilon)) \partial_y$.
\par We say that a limit cycle $\gamma_\varepsilon$ of system
(\ref{eq2}) {\em bifurcates from the origin} if it tends to the
origin (in the Hausdorff distance) as $\varepsilon \to 0$. We are
interested in giving a sharp upper bound for the number of limit
cycles which can bifurcate from the origin $p_0$ of system
(\ref{eq1}) under any analytic perturbation with a finite number
$p$ of parameters. The word sharp means that there exists a system
of the form (\ref{eq2}) with exactly that number of limit cycles
bifurcating from the origin, that is, the upper bound is
realizable. This sharp upper bound is called the {\em cyclicity}
of the origin $p_0$ of system (\ref{eq1}) and will be denoted by
${\rm Cycl}(\mathcal{X}_\varepsilon,p_0)$ all along this paper.
\par As defined in \cite{Andronov, HK, Takens}, a {\em Hopf bifurcation}, also
denoted by Poincar\'e-Andronov-Hopf bifurcation, is a bifurcation
in a neighborhood of a singular point like the origin of system
(\ref{eq1}). If the stability type of this point changes when
subjected to perturbations, then this change is usually
accompanied with either the appearance or disappearance of a small
amplitude periodic orbit encircling the equilibrium point.
\par
We remark that in all the examples considered in this paper we
take a perturbed system (\ref{eq2}) with the origin as singular
point, that is, $\bar{P}(0,0,\varepsilon)=
\bar{Q}(0,0,\varepsilon)\equiv 0$. Moreover, all the limit cycles
$\gamma_\varepsilon$ that bifurcate from the origin in our
examples surround the origin. \newline

We consider systems of the form (\ref{eq1}) where the origin $p_0$
is a focus singular point of the following three types:
non-degenerate, degenerate without characteristic directions and
nilpotent (the definitions are stated below). The study of the
cyclicity of a degenerate focus has been tackled in very few
sources; we mention the papers \cite{Chus-Ga, AnSaTs, CiGaMa,
DumRou91, GGLl06, Liu, Manosa} where some (partial) results can be
found. In relation with normal forms and integrability of
degenerate singular points, we cite \cite{ChGaGi, Liu, MatTei86}.
Our results do not establish that the cyclicity of this type of
singular points is finite but give an effective procedure to study
it. In the three mentioned types of focus points, we will consider
a change to (generalized) polar coordinates which embed the
neighborhood $\mathcal{U}$ of the origin into a cylinder $C\, = \,
\left\{ (r,\theta) \in \mathbb{R} \times \mathcal{S}^1 \, : \,
|r|<\delta \right\}$ for a certain sufficiently small value of
$\delta>0$. This change to polar coordinates is a diffeomorphism
in $\mathcal{U} - \{(0,0)\}$ and transforms the origin of
coordinates to the circle of equation $r=0$. In fact, the
neighborhood $\mathcal{U}$ is transformed into the half-cylinder
in which $r \geq 0$, but we can consider the extension to the
values in which $r<0$ by using several symmetries of the
considered (generalized) polar coordinates. In these new
coordinates, system (\ref{eq1}) can be seen as a differential
equation over the cylinder $C$ of the form:
\begin{equation}
\frac{dr}{d\theta} \, = \, \mathcal{F}(r,\theta), \label{eq3}
\end{equation} where $\mathcal{F}(r,\theta)$ is an analytic
function in $C$. The circle $r=0$ needs to be a particular
periodic solution of the equation (\ref{eq3}) and, therefore,
$\mathcal{F}(0,\theta) \equiv 0$ for all $\theta \in
\mathcal{S}^1$. \par Throughout the rest of the paper, we consider
an inverse integrating factor $V(r,\theta)$ of equation
(\ref{eq3}). We recall that an {\em inverse integrating factor} of
equation (\ref{eq3}) is a function $V: C \to \mathbb{R}$ of class
$\mathcal{C}^1(C)$, which is non locally null and which satisfies
the following partial differential equation:
\[ \frac{\partial V(r,\theta)}{\partial \theta} \, + \, \frac{\partial
V(r,\theta)}{\partial r} \, \mathcal{F}(r,\theta) \, = \,
\frac{\partial \mathcal{F}(r,\theta)}{\partial r} \, V(r,\theta).
\] We remark that since $V(r,\theta)$ is a function defined over
the cylinder $C$ it needs to be $T$--periodic in $\theta$, where
$T$ is the minimal positive period of the variable $\theta$, that
is, we consider the circle $\mathcal{S}^1 \, = \, \mathbb{R} /
[0,T]$. The function $V(r,\theta)$ is smooth
($\mathcal{C}^{\infty}$) and non--flat in $r$ in a neighborhood of
$r=0$. The existence of an inverse integrating factor
$V(r,\theta)$ with this regularity is proved in \cite{EncPer}
using the result of \cite{Yakovenko}, see also Lemma \ref{lemv2}
of the present paper. A characterization of when $V(r,\theta)$ is
analytic in a neighborhood of $r=0$ is given in \cite{EncPer}.

Let us consider the Taylor expansion of the function $V(r,\theta)$
around $r=0$: $ V(r,\theta) \, = \, v_m(\theta) \, r^m \, + \,
\mathcal{O}(r^{m+1}), $ where $v_m(\theta) \not\equiv 0$ for
$\theta \in \mathcal{S}^1$ and $m$ is an integer number with $m
\geq 0$. As we will see in the following section, in fact,
$v_m(\theta) \neq 0$ for all $\theta \in \mathcal{S}^1$, cf. Lemma
\ref{lemv1}. We say that $m$ is the {\em vanishing multiplicity}
of $V(r,\theta)$ on $r=0$. The aim of this paper is to show the
correspondence between this vanishing multiplicity $m$ and the
cyclicity ${\rm Cycl}(\mathcal{X}_\varepsilon,p_0)$ of the origin
$p_0$ of system (\ref{eq1}). \par One of our hypothesis is that
the origin of system (\ref{eq1}) is a focus, and thus, we obtain
that the circle $r=0$ is an isolated periodic orbit (i.e. a limit
cycle) of equation (\ref{eq3}). This hypothesis implies that there
can only exist one inverse integrating factor $V(r,\theta)$ smooth
and non--flat in $r$ in a neighborhood of $r=0$, up to a nonzero
multiplicative constant, see Lemma \ref{lemv2}. The uniqueness of
$V(r,\theta)$ implies that the number $m$ corresponding to the
vanishing multiplicity of $V(r,\theta)$ on $r=0$ is well--defined.
\newline

The existence of $V(r,\theta)$ gives the existence of an inverse
integrating factor $V_0(x,y)$ for system (\ref{eq1}) by undoing
the change to (generalized) polar coordinates. We recall that
$V_0: \mathcal{U} \to \mathbb{R}$ is said to be an {\em inverse
integrating factor} of system (\ref{eq1}) if it is of class
$\mathcal{C}^1(\mathcal{U})$, it is not locally null and it
satisfies the following partial differential equation:
\[ P(x,y) \, \frac{\partial V_0(x,y)}{\partial x} \, + \,
Q(x,y) \, \frac{\partial V_0(x,y)}{\partial y} \, = \, \left(
\frac{\partial P(x,y)}{\partial x} \, + \, \frac{\partial
Q(x,y)}{\partial y} \right) \, V_0(x,y). \] The function
$V_0(x,y)$, obtained from $V(r,\theta)$ by undoing the change to
polar coordinates, does not need to be smooth at the origin $p_0$
since the change to polar coordinates is a diffeomorphism except
in the origin. Besides, we are also interested in the problem of
the regularity of an inverse integrating factor $V_0(x,y)$ in a
neighborhood of the origin of system (\ref{eq1}), whenever it
exists, and we also analyze this question, see Theorem
\ref{thpoint} and Corollaries \ref{cornil} and \ref{corexist}.
\par The zero--set of $V_0(x,y)$, which we denote by
$V_0^{-1}(0):= \{ p\in \mathcal{U} \, : \,  V_0 (p) = 0\}$, is
formed by orbits of system (\ref{eq1}) and usually contains those
orbits which determine the dynamics of the system: singular
points, limit cycles and graphics, see \cite{BerroneGiacomini,
GaGiGr, Ga-Sh, GLV}. If there exists an inverse integrating factor
in a neighborhood of a limit cycle, then it is contained in
$V_0^{-1}(0)$, as it has been proved in \cite{GLV}. Since the
origin $p_0=(0,0)$ of (\ref{eq1}) is a focus, we have that
$V_0(0,0)=0$, as proved in \cite{BerroneGiacomini, Ga-Sh}.
Moreover, the set $V_0^{-1}(0)$ contains all the singular points
of (\ref{eq1}) having at least one parabolic or elliptic sector in
its domain of definition, as it is also proved in
\cite{BerroneGiacomini, Ga-Sh}.
\newline

In \cite{GaGiGr}, we proved that the vanishing multiplicity of an
analytic $V(r, \theta)$ on the limit cycle $r=0$ coincides with
the multiplicity (and, thus, the cyclicity) of $r=0$ as an orbit
of equation (\ref{eq3}). The statements and the proofs given in
\cite{GaGiGr} can be repeated verbatim with the weaker assumption
that $V(r,\theta)$ is smooth and non--flat in $r$ in a
neighborhood of $r=0$, see Corollary \ref{corvpi} in the following
section. The multiplicity of the limit cycle $r=0$ is related to
the cyclicity of the origin of system (\ref{eq1}) by the change
from polar coordinates. The aim of this work is to study the
cyclicity of the origin of system (\ref{eq1}) through the
vanishing multiplicity of $V(r, \theta)$ on $r=0$. We remark that
the multiplicity $m$ of $r=0$ as a limit cycle of equation
(\ref{eq3}) can also be established by successively solving the
variational equations. However, this method implies the
computation of iterated integrals of non--elementary periodic
functions. If we explicitly know an inverse integrating factor
$V_0(x,y)$ for system (\ref{eq1}), then we have an inverse
integrating factor $V(r,\theta)$ of (\ref{eq3}) and we can
immediately know the value of $m$ through the vanishing
multiplicity of $V(r, \theta)$ in $r=0$. On the other hand, since
the existence of an inverse integrating factor $V(r,\theta)$,
which is smooth and non--flat in $r$ in a neighborhood of $r=0$,
for equation (\ref{eq3}) is ensured (see Lemma \ref{lemv2}), we
have an alternative method to the variational equations to
determine the value of $m$.\newline

We are going to study the cyclicity of the origin of system
(\ref{eq1}) in the following three cases: the origin is a
non-degenerate focus, the origin is a degenerate focus without
characteristic directions and the origin is a nilpotent focus. The
two first cases are transformed to an equation of the form
(\ref{eq3}) using polar coordinates and the latter case using the
so-called generalized polar coordinates, see the definition below.
We state our results for each type of coordinates, separately. The
following Subsection \ref{sect11} contains the results related to
the non-degenerate and degenerate without characteristic
directions case. Subsection \ref{sect12} is devoted to the
nilpotent focus. In Subsection \ref{sect13} we analyze the
existence (and regularity conditions) of an inverse integrating
factor in a neighborhood of an isolated singular point of system
(\ref{eq1}). We give some general results on an equation
(\ref{eq3}) over a cylinder in Section \ref{sect2}. Finally,
Section \ref{sect3} contains the proofs of our main results.
\newline

\subsection{A focus without characteristic directions \label{sect11}}

We say that a focus at the origin of system (\ref{eq1}) is {\em
non-degenerate} if the linear part of system (\ref{eq1}) has
complex eigenvalues of the form $\alpha \, \pm \, \beta \, i$ with
$\alpha, \beta \in \mathbb{R}$ and $\beta \neq 0$. After a linear
change of coordinates and a rescaling of time, if necessary,
system (\ref{eq1}) can be written in the form:
\begin{equation} \label{eqnodeg}
\dot{x} \, = \, \zeta \, x \, -\, y \, + \, P_2(x,y), \qquad
\dot{y} \, = \, x \, + \, \zeta \, y \, + \, Q_2(x,y),
\end{equation} where $\zeta \in \mathbb{R}$ and $P_2(x,y)$ and $Q_2(x,y)$ are analytic
functions in a neighborhood of the origin without constant nor
linear terms.
\par
We say that the origin of system (\ref{eq1}) is a {\em degenerate
singular point} if the determinant associated to the linear part
of (\ref{eq1}) is zero. We consider a system (\ref{eq1}) of the
form:
\begin{equation} \label{eqdeg}
\dot{x} \, = \, p_d(x,y)\, +\, P_{d+1}(x,y), \qquad \dot{y} \, =
\, q_d(x,y)\, +\, Q_{d+1}(x,y),
\end{equation}
where $d  \geq 1$ is an odd number, $p_d(x,y)$ and $q_d(x,y)$ are
homogeneous polynomials of degree $d$ and $P_{d+1}(x,y),
Q_{d+1}(x,y) \in \mathcal{O}(\|(x,y)\|^{d+1})$, that is, they are
analytic functions in a neighborhood of the origin with order at
least $d+1$. We assume that $p_d^2(x,y)+q_d^2(x,y) \not\equiv 0$.
When $d > 1$, the origin of system (\ref{eqdeg}) is a degenerate
singular point.
\par A {\em characteristic direction} for the origin of system
(\ref{eqdeg}) is a linear factor in $\mathbb{R}[x,y]$ of the
homogeneous polynomial $xq_d(x,y)-yp_d(x,y)$. It is obvious that,
unless $xq_d(x,y)-yp_d(x,y) \equiv 0$, the number of
characteristic directions for the origin of system (\ref{eqdeg})
is less than or equal to $d+1$. If there are no characteristic
directions, then the origin is a monodromic singular point of
system (\ref{eqdeg}). We observe that the reciprocal is not true.
A singular point with characteristic directions can be monodromic.
The origin of system (\ref{eq1}) is a monodromic singular point if
it is either a center or a focus, see for instance \cite{Andronov,
GGG06, Manosa} for further information about monodromic singular
points and characteristic directions. We assume that the origin of
system (\ref{eqdeg}) is a focus without characteristic directions.
\newline

In the non-degenerate case ($d=1$), that the cyclicity of a focus
point is finite is well-known as well as several methods to
determine it. The most usual method is to compute the first
non-vanishing Liapunov constant and its order gives the cyclicity
of the focus. Indeed, the same method allows to study the limit
cycles which bifurcate from the origin in this case. When an
inverse integrating factor is known, we give a shortcut in the
study of the cyclicity as it can be given in terms of the
vanishing multiplicity of the inverse integrating factor at the
origin. \par As far as the authors know, the cyclicity of a
degenerate focus ($d>1$) of system (\ref{eqdeg}) is not proved to
be bounded. In fact, very few techniques appear in the literature
to tackle the cyclicity of this type of singular points. Usually,
polar coordinates are taken and the corresponding equation over
the cylinder is studied. Our approach is to take profit from the
knowledge of an inverse integrating factor to avoid the study of
the differential equation over the cylinder.
\newline

We remark that if $d=1$ and the origin of system (\ref{eqdeg}) is
a focus point (without characteristic directions), then it can be
written in the form (\ref{eqnodeg}).  \par We use polar
coordinates, $(x,y) \, = \, (r\, \cos\theta, \, r \, \sin\theta)$,
in order to transform a neighborhood of the origin into the
cylinder with period $T=2 \pi$, and system (\ref{eqdeg}) into an
ordinary differential equation of the form (\ref{eq3}). \par In
relation with system (\ref{eq2}), an analytic  perturbation field
$(\bar{P}(x,y,\varepsilon), \bar{Q}(x,y,\varepsilon))$ is said to
have subdegree $s$ if $(\bar{P}(x,y,\varepsilon), \bar{Q}(x, y,
\varepsilon)) = \mathcal{O}( \|(x,y)\|^s )$. In this case, we
denote by $\mathcal{X}_\varepsilon^{[s]}\, = \, (P(x,y) \, + \,
\bar{P}(x,y,\varepsilon))\partial_x \, + \, (Q(x,y) \, + \,
\bar{Q}(x,y,\varepsilon)) \partial_y$ the vector field associated
to such a perturbation.
\begin{theorem} \label{thdeg}
We assume that the origin $p_0$ of system {\rm (\ref{eqdeg})} is
monodromic and without characteristic directions. Let
$V(r,\theta)$ be an inverse integrating factor of the
corresponding equation {\rm (\ref{eq3})} which has a Laurent
expansion in a neighborhood of $r=0$ of the form $ V(r, \theta) \,
= \, v_{m}(\theta) \, r^m \, + \, \mathcal{O}(r^{m+1}), $ with
$v_m(\theta) \not\equiv 0$ and $m \in \mathbb{Z}$.
\begin{itemize}
\item[{\rm (i)}] If $m \leq 0$ or $m$ is even, then the origin of
system {\rm (\ref{eqdeg})} is a center. \item[{\rm (ii)}] If the
origin of system {\rm (\ref{eqdeg})} is a focus, then $m \geq 1$,
$m$ is an odd number and the cyclicity ${\rm
Cycl}(\mathcal{X}_\varepsilon,p_0)$ of the origin of system {\rm
(\ref{eqdeg})} satisfies ${\rm Cycl}(\mathcal{X}_\varepsilon,p_0)
\geq  (m+d)/2-1$. In this case $m$ is the vanishing multiplicity
of $V(r,\theta)$ on $r=0$.
\begin{itemize}
\item[{\rm (ii.1)}] If, moreover, the focus is non--degenerate $(d=1)$,
then the aforementioned lower bound is sharp, that is, ${\rm
Cycl}(\mathcal{X}_\varepsilon,p_0) =  (m-1)/2$.
\item[{\rm (ii.2)}] If only perturbations whose subdegree is greater than or equal to $
d$ are considered, then the maximum number of limit cycles which
bifurcate from the origin is $(m-1)/2$, that is, ${\rm
Cycl}(\mathcal{X}_\varepsilon^{[d]},p_0) =  (m-1)/2$.
\end{itemize}
\end{itemize}
\end{theorem}
The proof of this theorem is given in Section \ref{sect3}.

\begin{remark} \label{remdeg} As we will see in the proof of this theorem, if there exists an inverse integrating factor
$V_0(x,y)$ of system {\rm (\ref{eqdeg})} such that $V_0(r \cos
\theta, \, r  \sin \theta)/r^d$ has a Laurent expansion in a
neighborhood of $r=0$, then the exponents of the leading terms of
$V_0(r \cos \theta, \, r  \sin \theta)/r^d$ and $V(r,\theta)$
coincide. Thus, the vanishing multiplicity $m$ can be computed
without passing the system to polar coordinates.
\end{remark}

We provide several examples of application of Theorem \ref{thdeg}.
\newline

\begin{example} \label{ex1} {\rm The following system \begin{equation} \label{eqez}
\begin{array}{lll}
\dot{x} & = & \displaystyle -y \left( (2\mu+1) x^2 + y^2 \right)
\, + \, x^3 \left( \lambda_1 x^2 + \lambda_2 (x^2+y^2)\right),
\vspace{0.2cm}
\\ \dot{y} & = &  \displaystyle x \left( x^2 + (1-2\mu) y^2 \right) \, + \, x^2
y \left( \lambda_1 x^2 + \lambda_2 (x^2+y^2)\right), \end{array}
\end{equation}
where $\mu$, $\lambda_1$ and $\lambda_2$ are real parameters,
appears in \cite{ZoLl}, where it is shown that the origin is a
focus for a non semi-algebraic set of values of $(\mu, \lambda_1,
\lambda_2)$ and it is a center otherwise. \par
 We have that this system is written in the form
(\ref{eqdeg}) with $d=3$ and it has no characteristic directions
as $\ xq_3(x,y)-yp_3(x,y)\, = \, (x^2+y^2)^2$. Easy computations
show that the function \[ V_0(x,y) \, = \, e^{\frac{-2 \mu
x^2}{x^2+y^2}} \, (x^2+y^2)^3 \] is an inverse integrating factor
of system (\ref{eqez}) which satisfies the hypothesis of our
Theorem \ref{thdeg} with $V(r,\theta) \, = \,
e^{-2\mu\cos^2\theta} \, r^3$. We remark that $V_0(x,y)$ is not an
analytic function in a neighborhood of $r=0$ unless $\mu=0$. As a
consequence of Lemma \ref{lemv2}, which we prove in Section
\ref{sect3}, we deduce that when the origin of system (\ref{eqez})
is a focus, there are no analytic inverse integrating factors
$\bar{V}_0(x,y)$ defined in a neighborhood of the origin. If it
existed, its transformation to polar coordinates would produce an
analytic inverse integrating factor $\bar{V}(r,\theta)$ defined in
a neighborhood of $r=0$ and different from $V(r,\theta)$, up to a
multiplicative nonzero constant. In \cite{ZoLl} it is shown that
there exist values of $(\mu, \lambda_1, \lambda_2)$, either with
$\mu=0$ or with $\mu \neq 0$, for which the origin of system
(\ref{eqez}) is a focus. Applying Theorem \ref{thdeg}, we deduce
that whenever the origin of system (\ref{eqez}) is a focus, then
its cyclicity is greater than or equal to $2$. In the particular
case of a perturbation of subdegree greater or equal than $3$, the
maximum number of limit cycles that bifurcate from the origin is
$1$.}
\end{example}

\begin{example} \label{ex2} {\rm Let $k$, $s$ be integers such that $s \geq 2k
\geq 0$ and consider the following differential system:
\begin{equation} \label{exes}
\dot{x} \, = \, -y \, (x^2+y^2)^k \, + \, x \, R_{s}(x,y), \qquad
\dot{y} \, = \, x \, (x^2+y^2)^k \, + \, y \, R_{s}(x,y),
\end{equation} where $R_{s}(x,y)$ is a homogeneous polynomial of
degree $s$. The origin of this system is a monodromic singular
point since there are no characteristic directions. We take polar
coordinates $x \, = \, r \cos \theta$, $y \, = \, r \sin \theta$
and system (\ref{exes}) reads for \[ \dot{r} \, = \, r^{s+1} \,
R_{s}(\cos \theta, \sin \theta), \quad \dot{\theta} \, = \,
r^{2k}, \] where we have used that $R_{s}(r \cos \theta, r \sin
\theta)\, = \, r^{s} R_{s}(\cos \theta, \sin \theta)$ since
$R_{s}(x,y)$ is a homogeneous polynomial of degree $s$. Hence, we
see that the origin of system (\ref{exes}) is a focus if, and only
if, the following integral is different from zero:
\[ \mathcal{L} \, = \, \int_{0}^{2\pi} R_{s} \left( \cos \theta,
\sin \theta \right) \, d \theta \, \neq \, 0. \] We remark that
if $s$ is an odd number, then the origin of system (\ref{exes})
is a center. Easy computations show that
\[ V_0(x,y) \, = \, (x^2+y^2)^{s/2+1} \] is an inverse integrating
factor for system (\ref{exes}) and that $V(r,\theta) \, = \,
r^{s+1-2k}$ is an inverse integrating factor of the ordinary
differential equation on a cylinder corresponding to system
(\ref{exes}). Thus, applying Theorem \ref{thdeg}, we deduce that
the cyclicity ${\rm Cycl}(\mathcal{X}_\varepsilon,p_0)$ of the
origin of system (\ref{exes}), with $\mathcal{L} \neq 0$, is $\
{\rm Cycl}(\mathcal{X}_\varepsilon,p_0) \, \geq \, s/2$. \par We
observe that when $k=0$, we have that the origin of system
(\ref{exes}) is a non--degenerate monodromic singular point and
${\rm Cycl}(\mathcal{X}_\varepsilon,p_0)=s/2$. The center problem
for any even value of $s$ is determined only by one Lyapunov
constant, namely $\mathcal{L}$, whereas the cyclicity of the
origin (in case it is a focus) is given by $s/2$. We see, in this
way, that the center problem and the determination of the
cyclicity for a non--degenerate monodromic singular point are
strongly related but are not equivalent problems.} \end{example}

\begin{example} \label{ex3} {\rm In this example we show that no limit cycles
bifurcate from a focus of a homogeneous system of degree $d$,
under perturbations of subdegree $\geq d$. We consider a
homogeneous system:
\begin{equation} \label{eqhom} \dot{x} \, = \, p_d(x,y), \qquad \dot{y} \, = \, q_d(x,y), \end{equation}
where $p_d(x,y)$ and $q_d(x,y)$ are homogeneous polynomials of
degree $d$ with $d \geq 1$. A focus of system (\ref{eqhom}) has no
characteristic directions because any linear real factor of $y
p_d(x,y) - x q_d(x,y)$ gives an invariant straight line of system
(\ref{eqhom}) through the origin. Since $y p_d(x,y) - x q_d(x,y)$
has no linear real factors, we have that $d$ is odd. In polar
coordinates system (\ref{eqhom}) becomes
\[ \dot{r} \, = \,  r^d \,
R_d(\theta), \qquad \dot{\theta} \, = \, r^{d-1} \,
F_d(\theta), \]
where
\[
\begin{array}{lll}
R_d(\theta) &=& p_d(\cos \theta, \sin \theta) \cos \theta + q_d(\cos \theta,\sin \theta) \sin \theta , \\
F_d(\theta) &=& q_d(\cos \theta, \sin \theta) \cos \theta -
p_d(\cos \theta,\sin \theta) \sin \theta .
\end{array}
\]
The hypothesis that the origin is a focus implies that there are
no characteristic directions. The non-existence of characteristic
directions is equivalent to $F_d(\theta) \neq 0$ for $\theta \in
[0, 2\pi)$. We see that the origin of system (\ref{eqhom}) is a
focus if, and only if, \[ \mathcal{L} \, = \, \int_0^{2 \pi}
\frac{R_d(\theta)}{F_d(\theta)} \, d \theta \, \neq \, 0. \] \par
On the other hand, easy computations show that $V_0(x,y)\, = \, y
p_d(x,y) - x q_d(x,y)$ is an inverse integrating factor for system
(\ref{eqhom}), where we have used Euler's Theorem for homogeneous
polynomials. Indeed $V(r,\theta)=r$ is an inverse integrating
factor for the corresponding differential equation on the
cylinder. By Theorem \ref{thdeg} we have that, when the origin of
system (\ref{eqhom}) is a focus (i.e. $\mathcal{L} \neq 0$), no
limit cycles can bifurcate from it under perturbations of
subdegree $\geq d$. \par However, if we take perturbations of
lower degree, there can appear limit cycles which bifurcate from
the origin. By Theorem \ref{thdeg}, we have that at least
$(d-1)/2$ limit cycles bifurcate from the origin of system
(\ref{eqhom}). \par For instance, let us consider the following
system:
\begin{equation} \dot{x} = p_3(x,y) = (x-y)(x^2+y^2) \ , \ \dot{y} = q_3(x,y) = (x+y)(x^2+y^2),
\label{ejfd} \end{equation} which has an unstable degenerate focus
at the origin. The function $V_0(x,y)=(x^2+y^2)^2$ is an inverse
integrating factor for this system and, as we are under the
hypothesis of Theorem \ref{thdeg}, we deduce that $m=1$ and that
no limit cycle can bifurcate from system (\ref{ejfd}) under
perturbations of subdegree $\geq 3$. Let us consider the following
perturbed system:
\[ \dot{x} = (x-y)(x^2+y^2) \, - \, \varepsilon (x+y) \ , \ \dot{y} = (x+y)(x^2+y^2) \, +\,  \varepsilon
(x-y), \] which has the invariant algebraic curve
$x^2+y^2=\varepsilon$. For $\varepsilon>0$, this algebraic curve
is a hyperbolic limit cycle of the system which bifurcates from
the origin. The function
$V_\varepsilon(x,y)=(x^2+y^2)(x^2+y^2-\varepsilon)$ is an inverse
integrating factor of this system. } \end{example}

\subsection{A nilpotent focus \label{sect12}}

We say that the origin of system (\ref{eq1}) is a {\em nilpotent
singular point} if it is a degenerate singularity that can be
written as:
\begin{equation}\label{eqnil}
\dot{x} = y + P_2(x,y) \ , \ \dot{y} = Q_2(x,y) \ ,
\end{equation}
with $P_2(x,y)$ and $Q_2(x,y)$ analytic functions near the origin
without constant nor linear terms. The following theorem is due to
Andreev \cite{Andreev} and it solves the monodromy problem for the
origin of system (\ref{eqnil}), that is, it determines when the
origin is a monodromic singular point.
\begin{theorem} \label{thandreev} {\sc \cite{Andreev}} Let $y=F(x)$ be the solution of $y + P_2(x,y)
= 0$ passing through $(0, 0)$. Define the functions $f(x) = Q_2(x,
F(x)) = a x^\alpha + \cdots$ with $a \neq 0$ and $\alpha \geq 2$
and $\phi(x) = (\partial P_2/\partial x \, + \, \partial Q_2 /
\partial y)(x,F(x))$. We have that either
$\phi(x) = b x^\beta + \cdots$ with $b \neq 0$ and $\beta \geq 1$
or $\phi(x) \equiv 0$. Then, the origin of {\rm (\ref{eqnil})} is
monodromic if, and only if, $a < 0$, $\alpha = 2 n-1$ is an odd
integer and one of the following conditions holds:
\begin{itemize}
\item[{\rm (i)}] $\beta > n-1$.

\item[{\rm (ii)}] $\beta = n-1$ and $b^2+4 a n < 0$.

\item[{\rm (iii)}] $\phi(x) \equiv 0$.
\end{itemize}
\end{theorem}

\begin{definition} \label{defnil} We consider a system of the form {\rm (\ref{eqnil})} with the
origin as a monodromic singular point. We define its {\em Andreev
number} $n \geq 2$ as the corresponding integer value given in
Theorem {\em \ref{thandreev}}. \end{definition}

We consider system (\ref{eqnil}) and we assume that the origin is
a nilpotent monodromic singular point with Andreev number $n$.
Then, the change of variables
\begin{equation}\label{change1}
(x,y) \, \mapsto \, (x, y-F(x)),
\end{equation}
where $F(x)$ is defined in Theorem \ref{thandreev}, and the
scaling
\begin{equation}\label{change2}
(x,y) \, \mapsto \, (\xi \, x, -\xi \, y),
\end{equation}
with $\xi = (-1/a)^{1/(2-2n)}$, brings system (\ref{eqnil}) into
the following analytic form for monodromic nilpotent singularities
\begin{equation}\label{eqnil2}
\dot{x}\, = \, y\, (-1 + X_1(x,y)), \quad \dot{y} \, =\, f(x) + y
\, \phi(x) + y^2\, Y_0(x,y),
\end{equation}
where $X_1(0,0)=0$, $f(x) = x^{2n-1} + \cdots$ with $n \geq 2$ and
either $\phi(x) \equiv 0$ or $\phi(x) = b x^\beta + \cdots$ with
$\beta \geq n-1$. We remark that we have relabelled the functions
$f(x)$, $\phi(x)$ and the constant $b$ with respect to the ones
corresponding to system (\ref{eqnil}). We recall, cf. Theorem
\ref{thandreev}, that when $\beta=n-1$ we also have that
$b^2-4n<0$.
\newline

We are going to transform system (\ref{eqnil2}) to an equation
over a cylinder of the form (\ref{eq3}). The transformation
depends on the Andreev number $n$ and it is given through the {\it
generalized trigonometric functions} defined by Lyapunov
\cite{Liapunov} as the unique solution $x(\theta) = {\rm Cs} \,
\theta$ and $y(\theta) = {\rm Sn}\, \theta$ of the following
Cauchy problem
\begin{equation}\label{Cauchy}
\frac{d x}{d \theta} \, =\, -y , \ \frac{d y}{d \theta}\, =\,
x^{2n-1}, \qquad x(0)=1 , \, y(0)=0  .
\end{equation}
We observe that, in the particular case $n=1$, the previous
definition gives the classical trigonometric functions. \par We
introduce in $\mathbb{R}^2 \backslash \{(0,0)\}$ the change to
{\it generalized polar coordinates}, $\,  (x,y) \mapsto
(r,\theta)$, defined by  \begin{equation} \label{change} x\, =\, r
\, {\rm Cs}\, \theta , \qquad y\, =\, r^n \, {\rm Sn}\, \theta.
\end{equation} In relation with this change, we say that a
polynomial $R(x,y) \in \mathbb{C}[x,y]$ is a {\em
$(1,n)$--quasihomogeneous polynomial of weighted degree} $w$ if
the following identity is satisfied:
\[ R\left(\lambda \, x, \, \lambda^n \, y \right) \, = \, \lambda^{w} \, R(x,y), \]
for all $(x,y,\lambda) \in \mathbb{R}^3$. We observe that a
homogeneous polynomial of degree $w$ is, with this definition, a
$(1,1)$-quasihomogeneous polynomial of weighted degree $w$. \par
Since we are going to use some properties of the generalized
trigonometric functions and the relations satisfied by
$(1,n)$--quasihomogeneous polynomials of weighted degree $w$, we
summarize them up in the following proposition. The proof of each
of its statements can be found in \cite{Liapunov}.
\begin{proposition} \label{propgen} {\sc \cite{Liapunov}} \
We fix an integer $n \geq 1$ and we consider $({\rm Cs} \, \theta,
\, {\rm Sn}\, \theta)$ the solution of the Cauchy problem {\rm
(\ref{Cauchy})}. The following statements hold. \begin{itemize}
\item[{\rm (a)}] The functions ${\rm Cs}\, \theta$ and ${\rm Sn}\, \theta$ are
$T_n$--periodic with $ \displaystyle T_n \, =\,  2 \, \sqrt{\frac{\pi}{n}} \, \frac{\Gamma\left(
\frac{1}{2 n} \right)}{\Gamma\left( \frac{n+1}{2 n} \right)} \, ,
$ \\ where $\Gamma(\cdot)$ denotes the Euler Gamma function.
\item[{\rm (b)}]  $\, \displaystyle {\rm Cs}^{2n} \theta \, + \, n \, {\rm Sn}^2 \theta \, = \, 1$ (the fundamental relation).
\item[{\rm (c)}]  $ \, \displaystyle {\rm Cs} \, (-\theta)  \, = \,  {\rm
Cs} \, \theta , \ \  {\rm Sn} \, (-\theta)  \, = \, - {\rm
Sn}\, \theta , $ \\ $ \, \displaystyle {\rm Cs} \, (\theta+T_n/2)  \, = \, - {\rm
Cs} \, \theta , \ \  {\rm Sn} \, (\theta+T_n/2)  \, = \, - {\rm
Sn}\, \theta$.
\item[{\rm (d)}] Euler Theorem for quasihomogeneous polynomials: if $R(x,y)$ is a $(1,n)$--quasihomogeneous
polynomial of weighted degree $w$, then
\[ x \, \frac{\partial R(x,y)}{\partial x} \, + \, n \,y \, \frac{\partial R(x,y)}{\partial y} \, = \, w \, R(x,y) . \]
\item[{\rm (e)}]  $ \, \displaystyle {\rm Cs} \, \varphi  \, = \, -\, {\rm
Cs} \, \theta , \ {\rm Sn} \, \varphi  \, = \, (-1)^n \, {\rm
Sn}\, \theta , $ where $\varphi \, = \, (-1)^{n+1} \left(\theta \,
+ \,T_n/2 \right)$. If $\, R(x,y)$ is a $(1,n)$--quasihomogeneous polynomial of weighted degree $w$, then
\[ R\left({\rm Cs} \, \varphi, {\rm Sn} \, \varphi  \right) \, = \, (-1)^{w}
\, R\left({\rm Cs} \, \theta, {\rm Sn} \, \theta \right). \]  In
particular, \[ R(-1,0) \, = \, R \left({\rm Cs} \left(
T_n/2\right), {\rm Sn} \left( T_n/2 \right) \right) \, = \,
(-1)^{w} R \left({\rm Cs}\, 0, {\rm Sn}\, 0 \right) \, = \,
(-1)^{w} R \left(1, 0 \right). \]
\end{itemize} \end{proposition}

Analogously to the case of a degenerate focus without
characteristic directions, we can also provide the maximum number
of limit cycles which can bifurcate from a nilpotent focus when
only certain perturbations are taken into account. In this sense,
and in relation with system (\ref{eq2}), we consider the following
definition, which will be used in the following Theorem
\ref{thnil}.
\begin{definition} \label{defnilqh}
An analytic perturbation vector field $(\bar{P}(x,y,\varepsilon),
\bar{Q}(x,y,\varepsilon))$ is said to be {\em
$(1,n)$--quasihomogeneous of weighted subdegrees $(w_x,w_y)$} if
$\bar{P}(\lambda x,\lambda^n y,\varepsilon) \, = \,
\mathcal{O}(\lambda^{w_x})$ and $\bar{Q}(\lambda x,\lambda^n
y,\varepsilon) \, = \, \mathcal{O}(\lambda^{w_y})$. In this case,
we denote by $\mathcal{X}_\varepsilon^{[w_x,w_y]}\, = \, (P(x,y)
\, + \, \bar{P}(x,y,\varepsilon))\partial_x \, + \, (Q(x,y) \, +
\, \bar{Q}(x,y,\varepsilon)) \partial_y$ the vector field
associated to such a perturbation. \end{definition}

We remark that the perturbative functions
$\bar{P}(x,y,\varepsilon), \bar{Q}(x,y,\varepsilon)$ do not need
to be $(1,n)$-quasihomogeneous of a certain degree, they just need
to have a $(1,n)$-quasihomogeneous subdegree high enough.

The following theorem is one of the main results of this work. The
symbol $\lfloor x \rfloor$ denotes the integer part of $x$.

\begin{theorem} \label{thnil}
We assume that the origin of system {\rm (\ref{eqnil})} is
monodromic with Andreev number $n$. Let $V(r,\theta)$ be an
inverse integrating factor of the corresponding equation {\rm
(\ref{eq3})} which has a Laurent expansion in a neighborhood of
$r=0$ of the form $ V(r, \theta) \, = \, v_{m}(\theta) \, r^m \, +
\, \mathcal{O}(r^{m+1}), $ with $v_m(\theta) \not\equiv 0$ and $m
\in \mathbb{Z}$.
\begin{itemize}
\item[{\rm (i)}] If $m \leq 0$ or $m+n$ is odd, then the origin of
system {\rm (\ref{eqnil})} is a center. \item[{\rm (ii)}] If the
origin of system {\rm (\ref{eqnil})} is a focus, then $m \geq 1$,
$m+n$ is even and its cyclicity ${\rm
Cycl}(\mathcal{X}_\varepsilon,p_0)$ satisfies $\, {\rm
Cycl}(\mathcal{X}_\varepsilon,p_0) \, \geq \, (m+n)/2\, -\, 1 $.
In this case, $m$ is the vanishing multiplicity of $V(r,\theta)$
on $r=0$. \item[{\rm (iii)}] If the origin of system {\rm
(\ref{eqnil2})} is a focus and if only analytic perturbations of
$(1,n)$--quasihomogeneous weighted subdegrees $(w_x,w_y)$ with
$w_x \geq n$ and $w_y \geq 2n-1$ are taken into account, then the
maximum number of limit cycles which bifurcate from the origin is
$\lfloor (m-1)/2 \rfloor$, that is, ${\rm
Cycl}(\mathcal{X}_\varepsilon^{[n,2n-1]},p_0)\, = \, \lfloor
(m-1)/2 \rfloor$.
\end{itemize}
\end{theorem}

The proof of this theorem is given in Section \ref{sect3}.

\begin{remark} \label{remnil} The proof of this theorem shows that if there exists an inverse
integrating factor $V_0^{*}(x,y)$ of system {\rm (\ref{eqnil2})}
such that $V_0^{*}(r \, {\rm Cs}\, \theta, \, r^n \, {\rm Sn}\,
\theta)/r^{2n-1}$ has a Laurent expansion in a neighborhood of
$r=0$, then the exponents of the leading terms of $V_0^{*}(r \,
{\rm Cs}\, \theta, \, r^n  \, {\rm Sn}\, \theta)/r^{2n-1}$ and
$V(r,\theta)$ coincide. Therefore, the value of $m$ can be
determined without performing the transformation of the system to
generalized polar coordinates.
\end{remark}

The following corollary establishes a necessary condition for
system (\ref{eqnil}) to have an analytic inverse integrating
factor $V_0(x,y)$ defined in a neighborhood of the origin.

\begin{corollary} \label{cornil}
We assume that the origin of system {\rm (\ref{eqnil})} is a
nilpotent focus with Andreev number $n$, and that there exists an
inverse integrating factor $V_0(x,y)$ of {\rm (\ref{eqnil})} which
is analytic in a neighborhood of the origin. Then, $n$ is odd.
\end{corollary}
The proofs of Theorem \ref{thnil} and Corollary \ref{cornil} are
given in Section \ref{sect3}. Before the proofs of the main
results of this paper we provide several examples of application
of Theorem \ref{thnil}.
\newline

We would like to remark that the change to system (\ref{eqnil2})
is not always necessary to arrive at an equation over a cylinder,
using generalized polar coordinates. The following proposition
establishes sufficient conditions for an analytic system in
cartesian coordinates (nilpotent or not) to be transformed to an
equation over a cylinder by the change to generalized polar
coordinates. We remark that given any analytic function $P(x,y)$
in a neighborhood of the origin and any positive integer number
$\tilde{n}$, we can always develop $P(x,y)$ as a series of $(1,
\tilde{n})$--quasihomogeneous polynomials. That is, we can always
define $(1, \tilde{n})$--quasihomogeneous polynomials $p_i(x,y)$
of weighted degree $i$ such that the following identity is
satisfied $P(x,y) \, = \, \sum_{i \geq 0} p_i(x,y)$.
\begin{proposition} \label{propdirect}
Let $\tilde{n} \geq 2$ be an integer number and consider an
analytic system of the form: \begin{equation} \label{eqdirect}
\dot{x} \, = \, \sum_{i \geq a} p_i(x,y), \qquad \dot{y} \, = \,
\sum_{i \geq b} q_i(x,y), \end{equation} where $p_i(x,y)$ and
$q_i(x,y)$ are $(1, \tilde{n})$--quasihomogeneous polynomials of
weighted degree $i$ and $p_a(x,y)$, $q_b(x,y)$ are not identically
null. If $a-1=b-\tilde{n} \geq 0$ and
\[ {\rm Cs} \, \theta \, q_b\left({\rm Cs} \, \theta, {\rm Sn} \,
\theta \right) - \tilde{n} \, {\rm Sn} \, \theta \, p_a\left({\rm
Cs} \, \theta, {\rm Sn} \, \theta \right) \neq 0 \] for all
$\theta \in [0,T_{\tilde{n}}]$, then the origin of system {\rm
(\ref{eqdirect})} is monodromic. Moreover, the change of
coordinates $x \, = \, r \, {\rm Cs} \, \theta$, $y \, = \,
r^{\tilde{n}} \, {\rm Sn} \, \theta$ brings the system to an
equation of the form {\rm (\ref{eq3})} which is analytic in a
neighborhood of its periodic solution $r=0$.
\end{proposition}
{\em Proof of Proposition {\rm \ref{propdirect}}.} The change of
variables $(x,y) \mapsto (r, \theta)$ gives
\begin{eqnarray*}
\dot{r} &=& \frac{x^{2\tilde{n}-1} P(x,y) + y Q(x,y)}{r^{2\tilde{n}-1}} = \frac{r^{2\tilde{n}-1}
{\rm Cs}^{2\tilde{n}-1} \theta \sum_{i \geq a} r^i \tilde{p}_i(\theta) + r^{\tilde{n}} {\rm Sn}\,
\theta \sum_{i \geq b} r^i \tilde{q}_i(\theta)}{r^{2\tilde{n}-1}} \ , \\
\dot{\theta} & = & \frac{x Q(x,y) -\tilde{n} y
P(x,y)}{r^{\tilde{n}+1}} = \frac{r {\rm Cs} \, \theta \sum_{i \geq
b} r^i \tilde{q}_i(\theta) -\tilde{n} r^{\tilde{n}} {\rm Sn}\,
\theta \sum_{i \geq a} r^i \tilde{p}_i(\theta)}{r^{\tilde{n}+1}} \
,
\end{eqnarray*}
where $\tilde{p}_i(\theta) := p_i({\rm Cs} \theta, {\rm Sn}
\theta)$, $\tilde{q}_i(\theta) := q_i({\rm Cs} \theta, {\rm Sn}
\theta)$. When $a \geq 1$ and $b \geq \tilde{n}$, this system is
analytic at $r=0$. In this case, we have that
\begin{eqnarray*}
\dot{r} &=& {\rm Cs}^{2\tilde{n}-1} \theta \sum_{i \geq a} r^i
\tilde{p}_i(\theta)\, +\, r \, {\rm Sn}\,
\theta \sum_{i \geq b} r^{i-\tilde{n}} \tilde{q}_i(\theta) \ , \\
\dot{\theta} & = & {\rm Cs}\,  \theta \sum_{i \geq b}
r^{i-\tilde{n}} \tilde{q}_i(\theta) \, -\, \tilde{n}\, {\rm Sn}\,
\theta \sum_{i \geq a} r^{i-1} \tilde{p}_i(\theta) \ .
\end{eqnarray*}
Moreover, we get that
$$
\dot{\theta} = \left\{ \begin{array}{lll} -\tilde{n}\, {\rm Sn}\,
\theta \, \tilde{p}_a(\theta) r^{a-1} + \cdots & \mbox{if} & a-1 >
b-\tilde{n} \ , \vspace{0.2cm}
\\ {\rm Cs}\, \theta \, \tilde{q}_b(\theta)\,  r^{b-\tilde{n}} + \cdots &
\mbox{if} & a-1 < b-\tilde{n} \ , \vspace{0.2cm} \\ \left( {\rm
Cs}\, \theta \ \tilde{q}_b(\theta) \, - \, \tilde{n} \ {\rm Sn}\,
\theta \ \tilde{p}_a(\theta)\right)   r^{a-1 }+ \cdots & \mbox{if}
& a-1 = b-\tilde{n} \ ,
\end{array} \right.
$$
where the dots denote terms of higher order in $r$. In this way,
when $a-1 = b-\tilde{n}$ and ${\rm Cs}\theta \ \tilde{q}_b(\theta)
- \tilde{n} \ {\rm Sn}\theta \ \tilde{p}_a(\theta) \neq 0$ for all
$\theta \in [0, T_{\tilde{n}})$, we have that $\dot{\theta}$ is
different from zero for all $|r|$ small enough. \bbox
\newline

We want to remark that not every system with a nilpotent
singularity at the origin satisfies the hypothesis of Proposition
\ref{propdirect}. Moreover, the value of the Andreev number $n$ of
a nilpotent system satisfying Proposition \ref{propdirect} does
not need to coincide with the value of $\tilde{n}$ which appears
in the statement of this proposition.

\begin{example} \label{ex4} {\rm Let $m$ and $n$ be two positive integers such
that $n \geq 2$, $m \geq 1$ and $m+n$ is even. We consider the
following system:
\begin{equation} \label{exen}
\dot{x} \, = \, y \, +\,  x \, R(x,y), \quad \dot{y} \, = \,
-x^{2n-1}\, + \, n\, y \, R(x,y), \end{equation} where $R(x,y)$ is
a $(1,n)$-quasihomogeneous polynomial of weighted degree $m+n-2$.
The origin of system (\ref{exen}) is a nilpotent singularity and
it is a monodromic point applying Theorem \ref{thandreev}. We
observe that the Andreev number of system (\ref{exen}) is $n$. The
change to generalized polar coordinates (\ref{change}) brings the
system to the form:
\[ \dot{r} \,  = \, r^{m+n-1} \, R({\rm Cs} \,
\theta, \, {\rm Sn} \, \theta), \qquad \dot{\theta} \, = \,
r^{n-1}, \] where we have used the fundamental relation ${\rm
Cs}^{2n} \theta \, + \, n \, {\rm Sn}^2 \theta \, = \, 1$. Hence,
the condition for the origin of system (\ref{exen}) to be a focus
is that
\[ \mathcal{L} \, = \,  \, \int_{0}^{T_n} R\left( {\rm Cs} \,
\theta, \, {\rm Sn} \, \theta \right) \, d\theta  \, \neq \, 0, \]
where $T_n$ is the period defined in Proposition \ref{propgen}.
For instance, if we take $R(x,y)=x^{m+n-2}$ we always have that
$\mathcal{L} \neq 0$ because the integrand is positive. The
symmetry properties of the generalized trigonometric functions
${\rm Cs} \, \theta$, ${\rm Sn} \, \theta$ stated in Proposition
\ref{propgen} imply that in case $m+n$ is odd, the value of
$\mathcal{L}$ is zero. \par Easy computations show that
\[ V_0(x,y) \, = \, \left(x^{2n} \, + \, n \, y^2
\right)^{\frac{m-1}{2n}+1}  \] is an inverse integrating factor
for system (\ref{exen}), where we have used Euler Theorem for
$(1,n)$-quasihomogeneous polynomials, cf. Proposition
\ref{propgen}. It is clear that since $m$ and $n$ are arbitrary
positive integers with the restriction of being of the same
parity, the function $V_0(x,y)$ does not need to be analytic in a
neighborhood of the origin. We observe that $V_0(x,y)$ is analytic
in a neighborhood of the origin if, and only if, there exists an
integer $k \geq 0$ such that $m=2kn+1$, which is an odd integer.
We observe that:
\begin{itemize}
\item When the origin of system (\ref{exen}) is a focus, then $m$ and $n$ have the same parity.
\item When system (\ref{exen}) has an analytic inverse integrating factor $V_0(x,y)$ defined in a
neighborhood of the origin, then $m=2kn+1$, which is an odd
integer. Therefore, accordingly to Corollary \ref{cornil}, we
deduce that in this case and if the origin is a focus, then $n$
needs to be an odd integer.
\end{itemize}
The function $V(r,\theta)$ defined in the statement of Theorem
\ref{thnil} is $V(r, \theta)\, = \, r^m$ which implies that the
cyclicity ${\rm Cycl}(\mathcal{X}_\varepsilon,p_0)$ of the origin
of system (\ref{exen}), when it is a focus, satisfies ${\rm
Cycl}(\mathcal{X}_\varepsilon,p_0) \, \geq \, (m+n)/2-1 $. \par
Another way to see that the origin of system (\ref{exen}) is
always a center in case that $m+n$ is odd, is by using Theorem
\ref{thnil}. } \end{example}

\begin{example} \label{ex5} {\rm We fix an integer $n \geq 2$ and we consider the following planar differential system:
\begin{equation}
\dot{x} \, = \, y \, - \, \nu_1 \, x^n, \qquad \dot{y} \, = \, -
x^{2n-1} \, + \, \nu_2 \, x^{n-1} y, \label{exex} \end{equation}
where $\nu_1$ and $\nu_2$ are real parameters such that $\nu_1 \,
\nu_2  -  1  <  0$ and $(\nu_2+n\, \nu_1)^2-4n<0$. We remark that
the Andreev number of the nilpotent monodromic singular point at
the origin of system (\ref{exex}) is $n$. System (\ref{exex}) is
the most general $(1,n)$--quasihomogeneous planar polynomial
differential system of weighted degree $n$ with a nilpotent
monodromic singularity at the origin (we have used Andreev's
Theorem \ref{thandreev}). The change to generalized polar
coordinates (\ref{change}) transforms system (\ref{exex}) in:
\[ \begin{array}{lll} \dot{r} & = & r^n \, \, {\rm Cs}^{n-1} \theta
\left( - \nu_1 \, {\rm Cs}^{2n} \theta + \nu_2 \, {\rm Sn}^2 \theta \right), \vspace{0.2cm} \\
\dot{\theta} & = & - r^{n-1} \left( 1 - (\nu_2
+ \nu_1 \, n) \, {\rm Cs}^n \theta \, {\rm Sn} \, \theta \right). \end{array} \] Hence, we have that the
origin of system (\ref{exex}) is a focus if, and only if, the
following integral is different from zero: \[ \mathcal{L} \,
= \, \int_{0}^{T_n} \frac{{\rm Cs}^{n-1} \theta \left( - \nu_1 \,
{\rm Cs}^{2n} \theta + \nu_2 \, {\rm Sn}^2 \theta \right)}{ 1 - (\nu_2 + \nu_1 \, n) \, {\rm Cs}^n \theta \,
{\rm Sn} \, \theta } \  d \theta \, \neq
\, 0,
\] where $T_n$ is the minimal positive period of these generalized
trigonometric functions, see Proposition \ref{propgen}. We define
$z(\theta) \, := \, 1 - (\nu_2 + \nu_1 \, n) \, {\rm Cs}^n \theta
\, {\rm Sn} \, \theta $ and we have that $z(\theta)>0$ for all
$\theta \in \mathbb{R}$ because
\[ z(\theta) \, = \, \left( {\rm Cs}^{n} \theta \, - \, \frac{\nu_2 + \nu_1 \,
n}{2} \, {\rm Sn} \, \theta\right)^2 \, + \, \frac{1}{4} \,
\left(4n-(\nu_2+n\, \nu_1)^2\right){\rm Sn}^2 \theta, \] where we
have used the fundamental relation ${\rm Cs}^{2n} \theta \, + \, n
\, {\rm Sn}^2 \theta \, = \, 1$ and that, under our hypothesis,
$(\nu_2+n\, \nu_1)^2-4n<0$. \par We remark that, using the
symmetries of the generalized trigonometric functions stated in
Proposition \ref{propgen}, we deduce that if $n$ is even, then
$\mathcal{L}$ vanishes. Therefore, a necessary condition for the
origin of system (\ref{exex}) to be a focus is that $n$ must be
odd. Indeed, easy computations using the properties of the
generalized trigonometric functions, show that
\[ \frac{{\rm Cs}^{n-1} \theta \left( - \nu_1 \,
{\rm Cs}^{2n} \theta + \nu_2 \, {\rm Sn}^2 \theta \right)}{ 1 -
(\nu_2 + \nu_1 \, n) \, {\rm Cs}^n \theta \, {\rm Sn} \, \theta }
\, = \, \frac{\nu_2-n\nu_1}{2n} \, \frac{{\rm Cs}^{n-1}
\theta}{z(\theta)} \, + \, \frac{z'(\theta)}{2n\, z(\theta)} \, .
\]
Hence, we deduce that \[ \mathcal{L} \, = \,
\frac{\nu_2-n\nu_1}{2n} \, \int_{0}^{T_n} \frac{{\rm Cs}^{n-1}
\theta}{z(\theta)} \, d\theta . \] Under the assumptions $n >2$
and odd, $\nu_1 \, \nu_2 -  1  <  0$ and $(\nu_2+n\,
\nu_1)^2-4n<0$, we have that the integrand in the previous
expression ${\rm Cs}^{n-1} \theta/z(\theta) \geq 0$ for any value
of $\theta \in [0,T_n]$. Therefore, we have that, under these
assumptions, the origin of system (\ref{exex}) is a focus if, and
only if, $\nu_2 \neq n \nu_1$. \par The function $V_0(x,y) \, = \,
n y \dot{x} - x \dot{y} \, = \, x^{2n} \, - \, \left( \nu_2 +
\nu_1 \, n\right) x^n y \, + \, n y^2$ is an inverse integrating
factor of system (\ref{exex}), which is a polynomial and, thus, it
is analytic in the whole plane. We observe that the parity of $n$
agrees with Corollary \ref{cornil}. $ $From the given expression
of $V_0(x,y)$, we deduce that $V(r,\theta)=r$ is an inverse
integrating factor of the equation (\ref{eq3}) corresponding to
system (\ref{exex}). Applying Theorem \ref{thnil}, we get that
$m=1$. Thus, if we consider an analytic perturbation which is
$(1,n)$-quasihomogeneous of weighted subdegrees $(w_x,w_y)$, with
$w_x \geq n$ and $w_y \geq 2n-1$, then no limit cycles can
bifurcate from the origin. Indeed, we have that the cyclicity of
the origin of system (\ref{exex}) is at least $(n-1)/2$.}
\end{example}

\subsection{On the existence of an inverse integrating factor \label{sect13}}

The following result is a summary and a generalization of several
results on the existence of a smooth and non--flat inverse
integrating factor $V_0(x,y)$ in a neighborhood of an isolated
singular point, see \cite{flows, EncPer,GiaVia}.
\begin{theorem} \label{thpoint}
Let the origin be an isolated singular point of {\rm (\ref{eq1})}
and let $\lambda, \mu \in \mathbb{C}$ be the eigenvalues
associated to the linear part of {\rm (\ref{eq1})}. If $\lambda
\neq 0$, then there exists a smooth and non--flat inverse
integrating factor $V_0(x,y)$ in a neighborhood of the origin.
\end{theorem}

We can ensure the existence of an inverse integrating factor with
stronger regularity in some particular cases. We recall that a
singular point is said to be {\em strong} if the function
$\partial P/\partial x +\partial Q/\partial y$ is not zero on it
and it is said to be {\em weak} otherwise. \par By using a
translation, we can assume that the origin is the singular point
under consideration. We say that the origin is {\em analytically
linearizable} if there exists an analytic, near-identity change of
variables such that the transformed vector field is linear. We say
that the origin is {\em orbitally analytically linearizable} if
there exists an analytic, near-identity change of variables such
that the transformed vector field is a linear multiplied by a
scalar unit function.
\begin{corollary} \label{corpoint}
Let the origin be an isolated singular point of {\rm (\ref{eq1})}
and let $\lambda, \mu \in \mathbb{C}$ be the eigenvalues
associated to the linear part of {\rm (\ref{eq1})}. Then, the
following statements hold.
\begin{enumerate}[{\rm (i)}]
\item {\rm (Strong focus)} If $\lambda \, = \, \alpha + i \beta$
and $\mu \, = \, \alpha - i \beta$ with $\alpha, \beta \in
\mathbb{R}\backslash \{0\}$, then $V_0(x,y)$ is analytic and it is
unique up to a multiplicative constant. \item {\rm (Center)} If
$\lambda \, = \, i \beta$ and $\mu \, = \, - i \beta$ with $\beta
\in \mathbb{R}\backslash \{0\}$ and the origin is a center, then
$V_0(x,y)$ is analytic.

\item {\rm (Linearizable point)} If the origin is (orbitally) analytically
linearizable and $\lambda \neq 0$, then $V_0(x,y)$ is analytic.

\item {\rm (Node)} If $\ \lambda,  \mu \in \mathbb{R}$ and $\
\lambda  \, \mu>0$, then $V_0(x,y)$ is analytic and it is unique
up to a multiplicative constant.

\end{enumerate}
\end{corollary}

We include here the references where the previous statements have
been proved. All the proofs are based upon the same idea: take the
transformation to normal form, which is smooth in the considered
cases, see \cite{DuLlAr} for instance. There is an analytic
inverse integrating factor for the vector field in normal form
(usually polynomial) and, thus, undoing the transformation, the
obtained inverse integrating factor is smooth. In this way, if the
singular point is isolated and with $\lambda \neq 0$, then its
(orbitally) linearizability implies the existence of an analytic
inverse integrating factor as it has been shown in \cite{flows}.
If $\lambda \neq 0$, we have that the origin is either:
\begin{itemize}
\item[-] A strong focus, and the existence of an analytic inverse
integrating factor is given in \cite{GiaVia, flows, EncPer}. A
strong focus is a particular case of a linearizable point.
\item[-] A weak focus, that is, $\lambda \, = \, i \beta$ and $\mu
\, = \, - i \beta$ with $\beta \in \mathbb{R}\backslash \{0\}$ and
the origin is not a center. The existence of a smooth inverse
integrating factor is given in \cite{EncPer}, where a
characterization of the existence of an analytic inverse
integrating factor and an example of a weak focus without an
analytic inverse integrating factor defined in a neighborhood of
it are also given. In the following paragraph we include a sketch
of the proof of this fact for the sake of completeness.
\item[-] A non-degenerate center, and the existence of an analytic
inverse integrating factor is given by the normal form of a center
given by Poincar\'e. We recall that this normal form implies that
a non-degenerate center is orbitally analytically linearizable.
\item[-] A (hyperbolic) node, that is
$\lambda, \mu \in \mathbb{R}$ and $\lambda \, \mu >0$, and the
fact that the transformation to normal form is analytic is proved
in \cite{Arnold}. Thus, there is an analytic inverse integrating
factor defined in a neighborhood of it. Indeed, in a neighborhood
of a node there cannot exist a first integral, see \cite{flows},
so that the analytic inverse integrating factor is unique up to a
multiplicative constant.
\item[-] A (hyperbolic) saddle, that is $\lambda,  \mu \in
\mathbb{R}$ and $\lambda \, \mu <0$, and the existence of a smooth
inverse integrating factor is shown in \cite{EncPer}. \item[-] A
semi--hyperbolic point, that is $\mu=0$ but $\lambda \neq 0$. When
it is isolated, the existence of a smooth inverse integrating
factor is given in \cite{EncPer}.
\end{itemize}
For the sake of completeness and in relation with the topic of
this paper, we give the normal form in the case of a
non--degenerate weak focus, that is, we consider the analytic
system (\ref{eqnodeg})$_{\{ \zeta=0\} }$ and we suppose that its
origin is a focus. In \cite{Belitskii}, it is shown the existence
of a smooth and non--flat transformation that brings the
considered system to the Birkhoff normal form: \begin{equation}
\label{birk} \dot{x} \, = \, -y + S_1(r^2) x - S_2(r^2)y, \qquad
\dot{y} \, = \, x \, + \, S_2(r^2) x + S_1(r^2) y, \end{equation}
where $S_1$ and $S_2$ are formal series on $r^2=x^2+y^2$. We use
Borel's Theorem, see for instance \cite{Narasimhan}, to ensure the
existence of smooth functions representing $S_1$ and $S_2$.
\begin{theorem} {\sc [Borel's Theorem]}
For every point $p \in \mathbb{R}^n$ and for every formal series
in $n$ variables, there exists a $\mathcal{C}^{\infty}$ function
$f$ defined in a neighborhood of $p$ whose Taylor series at $p$ is
equal to the given formal series.\end{theorem} Therefore, we have
a smooth change of coordinates which brings system
(\ref{eqnodeg})$_{\{ \zeta=0\} }$ to system (\ref{birk}). Smooth
changes of coordinates in normal form theory usually come from an
order-to-order change and may induce flat terms in the normal
form. However, since system (\ref{eqnodeg})$_{\{ \zeta=0\} }$ is
analytic in a neighborhood of the origin, all of these flat terms
can be removed by a suitable smooth change of coordinates. On the
contrary, if system (\ref{eqnodeg})$_{\{ \zeta=0\} }$ were only
smooth and with an infinite codimension focus at the origin (that
is a focus with all its Liapunov constants equal to zero), then
the normal form (\ref{birk}) would be only formal because we
cannot ensure the removal of all the flat terms. For example, the
Birkhoff normal form of system
\[ \dot{x} \, = \, -y + x\, \exp\left(\frac{-1}{x^2+y^2}\right), \qquad
\dot{y} \, = \, x \, + \, y \, \exp\left(\frac{-1}{x^2+y^2}\right)
\] is $\dot{x} \, = \, -y$, $\dot{y}\, = \, x$ but its origin is a
focus. Therefore, the flat terms cannot be removed by any smooth
change. This is an example of a smooth system with a focus of
infinite codimension. It is known from Poincar\'e that an analytic
system (\ref{eqnodeg}) has no foci of infinite codimension. In
this work we only consider analytic differential systems.
\par We take the analytic system (\ref{eqnodeg})$_{\{ \zeta=0\} }$ and we
perform the smooth change of coordinates which brings it to the
Birkhoff normal form (\ref{birk}). It is well known, that the
origin of system (\ref{eqnodeg})$_{\{ \zeta=0\} }$ is a center if,
and only if, $S_1 \equiv 0$. Since we are in the case that the
origin is a focus, we have that $S_1$ is not identically null and
easy computations show that $V_0(x,y) \, = \, (x^2+y^2) \,
S_1(x^2+y^2)$ is an inverse integrating factor for the Birkhoff
normal form. By undoing the change to the original system, we
obtain a smooth and non--flat inverse integrating factor. \par For
a strong focus, there is a unique analytic inverse integrating
factor, as it has been shown in \cite{flows}. For a weak focus,
when there is an analytic inverse integrating factor, then it is
unique up to a multiplicative constant. However, there may exist
other inverse integrating factors with lower regularity, as it is
shown, for instance, in the forthcoming example with system
(\ref{ejbh}).
\newline

The following result is a consequence of the results given in
Theorems \ref{thdeg} and \ref{thnil} and it ensures the existence
of an inverse integrating factor, of class at least
$\mathcal{C}^1$, in a neighborhood of certain degenerate focus
points, namely for the origin of system (\ref{eqdeg}) without
characteristic directions and the origin of system (\ref{eqnil}).

\begin{corollary} \label{corexist}
There exists an inverse integrating factor $V_0(x,y)$, of class at
least $\mathcal{C}^1$, in a neighborhood of the following two
types of singular points: a degenerate focus without
characteristic directions and a nilpotent focus.
\end{corollary}

The proof of this corollary is given in Section \ref{sect3}.

\section{Ordinary differential equations over a cylinder \label{sect2}}

This section is devoted to several results related with ordinary
differential equations of the form (\ref{eq3}) defined over a
cylinder $C  \, = \, \left\{ (r, \theta) \in \mathbb{R} \times
\mathcal{S}^1 \, : \, |r|<\delta \right\}$ for a certain
$\delta>0$ sufficiently small. We denote by $T$ the minimal
positive period of the variable $\theta$, that is, we consider the
circle $\mathcal{S}^1 \, = \, \mathbb{R} / [0,T]$. Thus, we
consider an ordinary differential equation of the form
(\ref{eq3}):
\[ \frac{dr}{d\theta} \, = \,
\mathcal{F}(r,\theta), \] where $\mathcal{F}(r,\theta)$ is an
analytic function on the cylinder $C$ and $\mathcal{F}(0,\theta)
\equiv 0$. We have that, by assumption, the circle $r=0$ is a
periodic orbit of equation (\ref{eq3}).
\par
We assume that equation (\ref{eq3}) has an inverse integrating
factor $V(r,\theta)$ which is analytic in a neighborhood of $r=0$.
All the results remain true if we assume that $V(r,\theta)$ is a
smooth function in $C$ which is non--flat in $r$ in a neighborhood
of $r=0$. We remark that $V(r,\theta)$ is a function over the
cylinder $C$ and, thus, $V(r,\theta)$ is $T$--periodic in
$\theta$. We recall that we have defined the {\em vanishing
multiplicity} $m$ of $V(r,\theta)$ on $r=0$ as the value such that
\begin{equation} \label{eqvm} V(r, \theta) \, = \, v_m(\theta) \, r^m \, + \,
\mathcal{O} \left( r^{m+1} \right), \end{equation} with
$v_m(\theta) \not\equiv 0$. The following lemma is already stated
and proved in \cite{GaGiGr}. We include here a proof for the sake
of completeness.
\begin{lemma} \label{lemv1} {\sc \cite{GaGiGr}} \
If $V(r, \theta)$ is an inverse integrating factor of {\rm
(\ref{eq3})} which is smooth and non--flat in $r$ in a
neighborhood of $r=0$ and with vanishing multiplicity $m$ over
$r=0$, then the function $v_m(\theta)$ defined in {\rm
(\ref{eqvm})} satisfies that $v_m(\theta) \neq 0$ for $\theta \in
[0,T)$. \end{lemma} {\em Proof.}  By hypothesis
$\mathcal{F}(0,\theta) \equiv 0$ and we define the function
$F_1(\theta)$ as the one which satisfies
\[ \mathcal{F}(r, \theta) \, = \, F_1(\theta) \, r \, + \,
\mathcal{O} \left( r^2 \right). \] We note that $F_1(\theta)$ may
be identically null. We have that $V(r, \theta)$ satisfies the
following partial differential equation \[ \frac{\partial
V(r,\theta)}{\partial \theta} \, + \, \frac{\partial
V(r,\theta)}{\partial r} \, \mathcal{F}(r,\theta) \, = \,
\frac{\partial \mathcal{F}}{\partial r} \, V(r,\theta), \] and
equating the coefficients of the order $r^m$ in the previous
identity, we get that $ \, v_m'(\theta) \, = \, (1-m) \,
F_1(\theta) \, v_m(\theta)$. Since $v_m(\theta) \not\equiv 0$, let
$\theta_0 \in [0,T)$ be such that $v_m(\theta_0) \neq 0$. We
deduce that
\begin{equation} \label{eqvm0} v_m(\theta) \, = \, v_m(\theta_0)
\, \exp \left\{ (1-m) \, \int_{\theta_0}^{\theta} F_1(\sigma) \, d
\sigma \right\} \, . \end{equation} Therefore, we conclude that
$v_m(\theta) \neq 0$ for $\theta \in [0,T)$. \bbox \newline

As we will see below, $m$ coincides with the multiplicity of the
limit cycle $r=0$ of equation (\ref{eq3}). We remark that the
integral $\int_{0}^{T} F_1(\sigma) \, d \sigma $ is the
characteristic exponent of the periodic orbit $r=0$ in equation
(\ref{eq3}). In particular, either $m=1$ or $\int_{0}^{T}
F_1(\sigma) \, d \sigma \, = \, 0$. Thus, from the formula
(\ref{eqvm0}), we confirm that $v_m(\theta)$ is always a
$T$-periodic function. \newline

The following lemma establishes the existence and uniqueness of
$V(r, \theta)$ when the periodic orbit $r=0$ is isolated, that is,
when $r=0$ is a limit cycle of (\ref{eq3}). The existence is
proved in \cite{EncPer} and the uniqueness is also stated and
proved in \cite{GaGiGr}.
\begin{lemma} \label{lemv2} {\sc \cite{EncPer, GaGiGr}} \
If the circle $r=0$ is a limit cycle of {\rm (\ref{eq3})}, then
there exists an inverse integrating factor $V(r, \theta)$ of {\rm
(\ref{eq3})} which is smooth and non--flat in $r$ in a
neighborhood of $r=0$. Indeed, $V(r, \theta)$ is unique, up to a
nonzero multiplicative constant.
\end{lemma}
{\em Proof.} We recall here the main ideas to prove this
statement, for the sake of completeness. \par Let us assume that
$r=0$ is a limit cycle of multiplicity $m$ of equation
(\ref{eq3}). When $m = 1$, we say that $r=0$ is {\em hyperbolic}.
Following the ideas given in \cite{EncPer}, and by a result in
\cite{Yakovenko}, in a neighborhood of $r=0$, we can consider a
smooth, and non--flat in $\rho$ in a neighborhood of $\rho=0$,
change of coordinates $(r, \theta) \to (\rho, \tau)$ which takes
equation (\ref{eq3}) to
\[ \begin{array}{lll} \displaystyle \frac{d\rho}{d\tau} \, = \,
\lambda \, \rho \mbox{\quad with $\lambda \neq 0$}, & \ \mbox{if}
\ & m=1; \vspace{0.2cm} \\ \displaystyle \frac{d\rho}{d\tau} \, =
\, \rho^m \, + \, a \, \rho^{2m-1} \mbox{\quad with $a \in
\mathbb{R}$}, & \ \mbox{if} \ & m>1. \end{array} \] In the case
that $r=0$ is hyperbolic ($m=1$), we have that the change of
coordinates is, indeed, analytic in a neighborhood of $r=0$. We
remark that the function $\bar{V}(\rho,\tau)$ defined by:
\[ \bar{V}(\rho,\tau) \, : = \, \left\{ \begin{array}{lll} \displaystyle
\rho & \ \mbox{if} \ & m=1; \vspace{0.2cm} \\ \displaystyle \rho^m
\, + \, a \, \rho^{2m-1} & \ \mbox{if} \ & m>1; \end{array}
\right.
\] is an inverse integrating factor of the latter equation. By
undoing the change of coordinates, we have a smooth, and non--flat
in $r$ in a neighborhood of $r=0$, inverse integrating factor
$V(r,\theta)$ for equation (\ref{eq3}). In the case that $r=0$ is
hyperbolic, we have an analytic inverse integrating factor in a
neighborhood of $r=0$.
\par To show the uniqueness of $V(r,\theta)$, let us assume that there exist two linearly independent inverse
integrating factors $V(r, \theta)$ and $\tilde{V}(r, \theta)$ of
equation (\ref{eq3}) which are both smooth and non--flat in $r$ in
a neighborhood of $r=0$. We assume that $V(r, \theta)\, = \,
v_m(\theta) \, r^m \, + \, \mathcal{O}(r^{m+1})$ and $\tilde{V}(r,
\theta)\, = \, \tilde{v}_{\tilde{m}}(\theta) \, r^{\tilde{m}} \, +
\, \mathcal{O}(r^{\tilde{m}+1})$ and that $m \geq \tilde{m}$. We
have that the function on the cylinder $C$ defined by $H(r,
\theta) \, := \, V(r, \theta) / \tilde{V}(r, \theta)$ is not
locally constant, smooth in $r$ and of class $\mathcal{C}^1$ in
$\theta$, by using the Lemma \ref{lemv1}. \par If $m>\tilde{m}$,
we have that $H(r,\theta)$ is constant equal to $0$ all over the
circle $r=0$ and if $m \, = \, \tilde{m}$, using the proof of
Lemma \ref{lemv1}, we have that $H(0,\theta) \, = \,
v_m(0)/\tilde{v}_m(0)$, from which we deduce that $H(r,\theta)$
takes a constant value all over the circle $r=0$. Moreover, this
function $H(r, \theta)$ is a first integral of equation
(\ref{eq3}), since it satisfies that \[ \frac{\partial H(r,
\theta)}{\partial \theta} \, + \, \frac{\partial H(r,
\theta)}{\partial r} \, \mathcal{F}(r,\theta) \, = \, 0.
\] Thus, $H(r, \theta)$ is constant on each orbit of equation
(\ref{eq3}). If $r=0$ is a limit cycle, then the orbits in a
neighborhood of this circle accumulate on it. By continuity of
$H(r,\theta)$, this fact implies that $H(r, \theta)$ needs to take
the same value on any point in a neighborhood of $r=0$ in
contradiction with the fact that $H(r, \theta)$ is not locally
constant. \bbox \newline

The following example, which is described in page 219 of
\cite{BerroneGiacomini}, shows that the conditions for $V(r,
\theta)$ to be smooth and non--flat in $r$ in a neighborhood of
$r=0$ and $2 \pi$-periodic in $\theta$ are essential to have a
unique inverse integrating factor. We consider the planar
differential system \begin{equation} \dot{x} \, = \, -y +
x(x^2+y^2), \quad \dot{y} \, = \, x + y(x^2+y^2), \label{ejbh}
\end{equation} which has a non-degenerate (unstable) focus at the
origin. The following functions are two inverse integrating
factors of the system which are of class $\mathcal{C}^1$ in a
neighborhood of the origin \[ V_0(x,y)=(x^2+y^2)^2 \ \mbox{ and }
\ \bar{V_0}(x,y)=(x^2+y^2)^2
\sin\left(2\arctan\left(\frac{y}{x}\right)+\,
\frac{1}{x^2+y^2}\right).
\] In polar coordinates, this system reads for
$\dot{r} \, = \, r^3$, $\dot{\theta} \, = \, 1$, which has the
following inverse integrating factor $V_1(r,\theta)\, = \,
V_0(r\cos \theta, r \sin \theta)/r \, = \, r^3$ analytic in a
neighborhood of $r=0$. Moreover, the function $V_2(r,\theta) \, =
\, \bar{V_0}(r\cos \theta, r \sin \theta)/r \, = \, r^3 \, \sin(2
\theta + r^{-2})$ is an inverse integrating factor of class
$\mathcal{C}^1$ in $r \neq 0$.  Indeed, $V_3(r,\theta) \, = \, r +
2r^3 \theta$ is another inverse integrating factor of the system
in polar coordinates, which is analytic in $r$ but not
$2\pi$--periodic in $\theta$.
\newline

We will show that when the periodic orbit $r=0$ is a limit cycle
of equation {\rm (\ref{eq3})}, the vanishing multiplicity of an
inverse integrating factor $V(r,\theta)$ is strictly positive.
\begin{lemma} \label{lemv0}
Let us consider a differential equation of the form {\rm
(\ref{eq3})} over a cylinder $C  \, = \, \left\{ (r, \theta) \in
\mathbb{R} \times \mathcal{S}^1 \, : \, |r|<\delta \right\}$ for a
certain $\delta>0$ and let us assume that $r=0$ is a periodic
orbit of the equation. We assume that $V(r,\theta)$ is an inverse
integrating factor for equation {\rm (\ref{eq3})} defined in
$C\backslash\{ r=0\}$ and which has a Laurent series of the form:
\[ V(r, \theta) \, = \, v_{m}(\theta) \, r^{m} \, + \, \mathcal{O}\left(r^{m+1}\right), \]
with $v_{m}(\theta) \not\equiv 0$ and $m \in \mathbb{Z}$. If $m
\leq 0$, the periodic orbit $r=0$ has a neighborhood filled with
periodic orbits, that is, it is not a limit cycle. \end{lemma}
{\em Proof.} Let $\omega \, = \, dr \, - \,
\mathcal{F}(r,\theta)\, d\theta$ be the Pfaffian $1$--form
associated to equation (\ref{eq3}). Since $V(r,\theta)$ is an
inverse integrating factor of the equation (\ref{eq3}), we have
that $\omega/V$ is a closed 1--form. We observe that the proof of
Lemma \ref{lemv1} also applies and, therefore, we have that
$v_{m}(\theta) \neq 0$ for any $\theta \in [0,T)$. In case that
$m\leq 0$, we have that $\omega/V$ is well-defined in the whole
cylinder $C$. Let us consider any non--contractible cycle in this
cylinder, for instance the cycle $r=0$. By virtue of De Rham's
Theorem, see \cite{Fulton}, we have that the closed $1$--form
$\omega/V$ is exact if, and only if, the value of the line
integral $ \int_{r=0} \omega/V$ is zero. This value is zero since
the oval $r=0$ is an orbit of the 1-form $\omega$ and, thus,
$\omega_{\left|r=0\right.} \, \equiv \, 0$. Hence, we have that
$\omega/V\, =\, d H$ for a certain $\mathcal{C}^2$ function
$H(r,\theta)$, which turns out to be a first integral of equation
(\ref{eq3}). The existence of this first integral implies that the
cycle $r=0$ is surrounded by periodic orbits, formed by the level
curves of $H$. \bbox \newline

In order to relate the vanishing multiplicity $m$ of $V(r,\theta)$
on $r=0$, which we have proved to be positive, and the cyclicity
of the focus at the origin of system (\ref{eq1}), we use a
previous result which is already stated and proved in
\cite{GaGiGr}. Our result gives an ordinary differential equation
for the Poincar\'e map associated to equation (\ref{eq3}) in terms
of the inverse integrating factor $V(r,\theta)$. If the minimal
positive period of $\theta$ in (\ref{eq3}) is denoted by $T$ and
$\Psi(\theta;r_0)$ is the flow of (\ref{eq3}) with initial
condition $\Psi(0;r_0)=r_0$, we recall that the Poincar\'e map
$\Pi: \Sigma \subseteq \mathbb{R} \to \mathbb{R}$ is defined as
$\Pi(r_0) \, = \, \Psi(T;r_0)$. We have that $\Pi$ is an analytic
diffeomorphism defined in a neighborhood $\Sigma$ of $r_0=0$ whose
fixed points correspond to periodic orbits of the equation
(\ref{eq3}). Since we have that $r=0$ is a limit cycle of equation
(\ref{eq3}), we deduce that the Poincar\'e map is not the identity
map. We recall that the limit cycle $r=0$ is said to have {\em
multiplicity $k$} if the expansion of the analytic Poincar\'e map
in a neighborhood of $r_0=0$ is of the form $ \Pi(r_0) \, = \, r_0
\, + \, c_k \, r_0^k \, + \, \mathcal{O}(r_0^{k+1}), $ where $c_k
\neq 0$. The multiplicity of the limit cycle $r=0$ in equation
(\ref{eq3}) allows us to determine the cyclicity of the focus at
the origin of system (\ref{eq1}), as we will see below. \par We
state the result proved in \cite{GaGiGr} in the terms used in the
present paper.
\begin{theorem} {\sc \cite{GaGiGr}} \label{thvpi} \ Let us assume that $V(r,\theta)$ is an inverse integrating factor
of equation {\rm (\ref{eq3})} defined in a neighborhood of the
periodic orbit $r=0$, whose minimal positive period is denoted by
$T$. We consider $\Pi(r_0)$ the Poincar\'e map associated to the
periodic orbit $r=0$ of equation {\rm (\ref{eq3})}. Then, the
following identity holds:
\begin{equation}
V \left(\Pi(r_0),T\right) \, = \, V \left(r_0,0\right) \,
\Pi'(r_0). \label{eqvpi}
\end{equation}
\end{theorem}

As a consequence of equation (\ref{eqvpi}), we can prove the
following result.
\begin{corollary} {\sc \cite{GaGiGr}} \label{corvpi} \ Let us assume that $V(r,\theta)$ is an inverse integrating factor
of equation {\rm (\ref{eq3})} which is smooth and non--flat in $r$
in a neighborhood of the limit cycle $r=0$ and whose vanishing
multiplicity over it is $m$. Then, $r=0$ is a limit cycle of
multiplicity $m$.
\end{corollary}
{\em Proof.} Since $V(r,\theta)$ is assumed to be a function on
the cylinder $C$, we have that it needs to be $T$ periodic in
$\theta$. We have that $V(r_0,T) \, = \, V(r_0,0)$ and we consider
its development in a neighborhood of $r_0=0$: $V(r_0,0) \, = \,
\nu_m \, r_0^m \, + \, \mathcal{O}\left( r_0^{m+1} \right)$, where
$\nu_m \neq 0$. We observe that the index $m$ appearing in this
decomposition coincides with the vanishing multiplicity of
$V(r,\theta)$ in $r=0$ by Lemma \ref{lemv1}. Recalling that $ \Pi(r_0) \, = \, r_0 \, +
\, c_k \, r_0^k \, + \, \mathcal{O}(r_0^{k+1}), $ where $c_k \neq
0$, we consider equation (\ref{eqvpi}) and we
subtract $V\left(r_0,0\right)$ from both members to obtain that:
\[ \nu_m \left( \Pi(r_0)^m  - r_0^m \right) \, + \,
\mathcal{O}\left(r_0^{m+1}\right) \, = \, \left(\nu_m \, r_0^m \,
+ \, \mathcal{O}\left(r_0^{m+1}\right) \right) \left(k \, c_k \,
r_0^{k-1} \, + \, O(r_0^k) \right).
\]
The lowest order terms in both sides of the previous identity
correspond to $r_0^{m+k-1}$ and the equation of their coefficients
is $m\, c_k \, \nu_m = \, k \, c_k \, \nu_m$, which implies that
$k=m$. \bbox

\section{Proofs of the results \label{sect3}}

\subsection*{Proof of Theorem {\rm \ref{thdeg}}.}

The following lemma establishes the first step in the proof of Theorem \ref{thdeg}. We show that the transformation to polar coordinates $x=r \cos \theta$, $y=r \sin \theta$, of system (\ref{eqdeg}) gives an equation over a cylinder of the form {\rm
(\ref{eq3})}.

\begin{lemma} \label{lemdeg1}
We consider system {\rm (\ref{eqdeg})} with $d \geq 1$ and odd
integer and we assume there are no characteristic directions. Then
the transformation to polar coordinates brings system {\rm (\ref{eqdeg})} to an ordinary
differential equation over a cylinder.
\end{lemma}
{\em Proof.} In polar coordinates system (\ref{eqdeg}) becomes
\begin{equation} \label{eq5a}
\begin{array}{lllll}
\dot{r} & = & r^d \, R(r,\theta) & = & r^d \,
(R_d(\theta)+\mathcal{O}(r)), \vspace{0.2cm} \\
\dot{\theta} & = & r^{d-1} F(r,\theta) & = & r^{d-1} \,
(F_d(\theta) + \mathcal{O}(r)) ,
\end{array}
\end{equation}
where
\[
\begin{array}{lll}
R_d(\theta) &=& p_d(\cos \theta, \sin \theta) \cos \theta + q_d(\cos \theta,\sin \theta) \sin \theta , \\
F_d(\theta) &=& q_d(\cos \theta, \sin \theta) \cos \theta -
p_d(\cos \theta,\sin \theta) \sin \theta .
\end{array}
\]
The hypothesis that there are no characteristic directions is
equivalent to say that $F_d(\theta) \neq 0$ for $\theta \in [0,
2\pi)$. We can, therefore, consider the ordinary differential
equation associated to the orbits of system (\ref{eq5a}) which
takes the form (\ref{eq3}):
\begin{equation} \label{eq5b} \frac{dr}{d\theta} \, = \, \frac{r
\, R(r,\theta)}{F(r,\theta)} \, . \end{equation} \bbox \newline

The following lemma establishes that the center problem for the
origin of system (\ref{eqdeg}) is equivalent to the determine when
the circle $r=0$ is contained in a period annulus for equation
(\ref{eq5b}).

\begin{lemma} \label{lemdeg1b}
We consider system {\rm (\ref{eqdeg})} with $d \geq 1$ and odd
integer and we assume there are no characteristic directions. The
circle $r=0$ is a limit cycle of equation {\rm (\ref{eq5b})} if,
and only if, the origin of system {\rm (\ref{eqdeg})} is a focus.
\end{lemma}
{\em Proof.} The origin of system (\ref{eqdeg}) is transformed to
the periodic orbit $r=0$ of equation (\ref{eq5b}) by the
transformation to polar coordinates. This transformation gives a
one-to-one correspondence between any point in a punctured
neighborhood of the origin in the plane $(x,y)$ and a cylinder $\{
(r,\theta) \, : \, 0<r<\delta, \theta \in \mathcal{S}^1\}$ for
$\delta>0$ sufficiently small. Thus, any orbit spiraling from or
towards the origin of system (\ref{eqdeg}) is transformed to an
orbit spiraling (from or towards) the circle $r=0$ in equation
(\ref{eq5b}). \bbox \newline

We have a symmetry for equation (\ref{eq5b}) which is inherited by
the symmetries of the polar coordinates.
\begin{lemma} \label{lemsym}
Let us consider a planar $\, \mathcal{C}^1$ differential system
$\dot{x} \, = \, P(x,y)$, $\dot{y} \, = \, Q(x,y)$, and perform
the change to polar coordinates $x \, = \, r \, \cos \theta$, $y
\, = \, r\, \sin \theta$. The resulting system $\dot{r} \, = \,
\Xi(r, \theta)$, $\dot{\theta} \, = \, \Theta(r, \theta)$ is
invariant under the change of variables $ (r, \theta) \, \mapsto
\, \left(-r, \theta + \pi\right). $
\end{lemma}
{\it Proof.} We observe that the monomials $r  \cos \theta$ and $r
\sin \theta$ are invariant by the change of variables. Since
\begin{eqnarray*} \dot{r} & = & \displaystyle \frac{r  \cos \theta \, P(r  \cos \theta,r \sin
\theta) \, +\,  r \sin \theta \, Q(r  \cos \theta,r \sin
\theta)}{r} , \vspace{0.2cm}
\\
\dot{\theta} & = & \frac{r  \cos \theta \, Q(r  \cos \theta,r \sin
\theta)\,  - \, r \sin \theta \, P(r  \cos \theta,r \sin
\theta)}{r^2} ,
\end{eqnarray*}
the result follows.  \bbox

\begin{lemma} \label{lemddeg}
We assume that the origin of system {\rm (\ref{eqdeg})} is a focus
without characteristic direction and we consider the corresponding
equation {\rm (\ref{eq5b})} with associated Poincar\'e map
$\Pi(r_0) \, = \, r_0 \, + \, c_m r_0^m \, + \,
\mathcal{O}(r_0^{m+1}),$ with $c_m \neq 0$. Then, $m$ is odd.
 \end{lemma}
{\em Proof.} By Lemma \ref{lemdeg1b}, we have that the origin of
system (\ref{eqdeg}) is a focus if, and only if, the circle $r=0$
is a limit cycle of equation {\rm (\ref{eq5b})}. We assume that
the circle $r=0$ is a limit cycle with multiplicity $m$ and we
consider its associated Poincar\'e map $\Pi(r_0)$. By Lemma
\ref{lemsym} we have that equation (\ref{eq5b}) has the discrete
symmetry $(r, \theta) \mapsto (-r, \theta + \pi)$ because it comes
from a system in cartesian coordinates (\ref{eqdeg}). This
symmetry implies that $r=0$ is either a stable or an unstable
limit cycle. Thus, $m$ needs to be odd.\bbox \newline

The next lemma states that, under our hypothesis, we have an
inverse integrating factor $V(r, \theta)$ for equation
(\ref{eq5b}) which is smooth and non--flat in $r$ in a
neighborhood of $r=0$ and such that $V(0, \theta) \equiv 0$.

\begin{lemma} \label{lemdeg2} We consider system {\rm (\ref{eqdeg})} with $d \geq 1$ and odd
integer and we assume there are no characteristic directions.
\begin{itemize}
\item[{\rm (i)}] If system {\rm (\ref{eqdeg})} has an inverse integrating
factor $V_0(x,y)$ defined in a neighborhood of the origin, then
the function defined by \[ V(r, \theta) \, := \, \frac{V_0(r \cos
\theta, r \sin \theta)}{r^d\, F(r,\theta)} \, \] is an inverse
integrating factor for equation {\rm (\ref{eq5b})} in $r \, \neq
\, 0$.
\item[{\rm (ii)}] Let $V(r,\theta)$ be an inverse integrating factor of equation {\rm (\ref{eq5b})} which has a Laurent expansion in a neighborhood of $r=0$ of the form $ V(r, \theta) \, = \, v_{m}(\theta) \, r^m \, + \, \mathcal{O}(r^{m+1}), $ with $v_m(\theta) \not\equiv 0$ and $m \in \mathbb{Z}$. Then, if $m \leq 0$ the origin of system {\rm
(\ref{eqdeg})} is a center.
\end{itemize}
\end{lemma} {\em Proof.} {\rm (i)} Since the jacobian to polar coordinates
is $r$, we have that the function $V_0(r \cos \theta, r \sin
\theta)/r$ is an inverse integrating factor for system
(\ref{eq5a}). We see that equation (\ref{eq5b}) is the equation of
the orbits associated to system (\ref{eq5a}) and, therefore, the
function $V(r, \theta) \, := \,V_0(r \cos \theta, r \sin
\theta)/(r^d\, F(r,\theta))$ is an inverse integrating factor of
(\ref{eq5b}). We remark that $V(r,\theta)$ does not need to be
well-defined in a neighborhood of $r=0$. Thus, we only have, at
the moment, that it is an inverse integrating factor for equation
{\rm (\ref{eq5b})} in $r \neq 0$.
\par {\rm (ii)} If $m\leq 0$, we are under the hypothesis  of Lemma \ref{lemv0} and we conclude that the cycle
$r=0$ is surrounded by periodic orbits, which give rise to a
neighborhood of the origin of system (\ref{eqdeg}) filled with
periodic orbits. Therefore, the origin of system {\rm
(\ref{eqdeg}) is a center. \bbox
\newline

The previous Lemma \ref{lemdeg2}, together with the statements
given in Lemmas \ref{lemdeg1} and \ref{lemdeg1b}, establishes that
under the hypothesis of Theorem \ref{thdeg}, we have that $m \geq
1$. Moreover, by Lemma \ref{lemddeg} and Corollary \ref{corvpi},
we have that $m$ needs to be odd. We have proved statement (i) in
Theorem \ref{thdeg}. \newline

Assuming that the origin of system (\ref{eqdeg}) is a focus, the
following step of the proof is to relate the multiplicity of the
limit cycle $r=0$ of equation (\ref{eq5b}) and the cyclicity of
the origin of system (\ref{eqdeg}).

\begin{lemma} We consider system {\rm (\ref{eqdeg})} with $d \geq 1$ and odd
integer, we assume there are no characteristic directions and that
the origin $p_0$ is a focus with cyclicity ${\rm
Cycl}(\mathcal{X}_\varepsilon,p_0)$. We consider the corresponding
equation {\rm (\ref{eq5b})} and we assume that $r=0$ is a limit
cycle of multiplicity $m=2k+1$. Then, ${\rm
Cycl}(\mathcal{X}_\varepsilon,p_0) \geq (m+d)/2-1$. Moreover:
\begin{itemize}
\item[{\rm 1.}] When $d=1$, we have ${\rm
Cycl}(\mathcal{X}_\varepsilon,p_0) \, = \, k$. \item[{\rm 2.}] If
we only consider perturbations of system {\rm (\ref{eqdeg})} whose
subdegree is $\geq d$, then at most $k$ limit cycles can bifurcate
from the origin of system {\rm (\ref{eqdeg})}, that is, ${\rm
Cycl}(\mathcal{X}_\varepsilon^{[d]},p_0) \, = \, k$.
\end{itemize}
\end{lemma}
{\em Proof.} The particular case $d=1$ of statement {\rm 1.} and
system (\ref{eqdeg}) with a focus at the origin is proved in
Theorem 40, Ch. IX in \cite{Andronov} (page 254), see also
\cite{Takens}.
\par We first provide an example of a perturbation of equation
(\ref{eq5b}), with $m$ limit cycles bifurcating from $r=0$, whose
transformation to cartesian coordinates gives a perturbation of
system (\ref{eqdeg}) of the form (\ref{eq2}) with exactly
$(m-1)/2=k$ limit cycles bifurcating from the origin. This example
shows that ${\rm Cycl}(\mathcal{X}_\varepsilon,p_0) \geq k$. \par
We take the corresponding equation (\ref{eq5b}) and we assume that
$r=0$ is a limit cycle of multiplicity $m$, which is an odd
integer with $m \geq 1$. We consider the associated system
(\ref{eq5a}) from which (\ref{eq5b}) comes from, we define
$k=(m-1)/2$ and we perturb system (\ref{eq5a}) in the following
way:
\begin{equation} \dot{r} \, = \, r^d \, R(r, \theta)
\, + \, \sum_{i=0}^{k-1} \varepsilon^{k-i} \, a_i \, r^{2i+d},
\qquad \dot{\theta} \, = \, r^{d-1} \, F(r, \theta), \label{eq5p}
\end{equation} with the convention that if $k=0$ no
perturbation term is taken. The real constant $\varepsilon$ is the
perturbation parameter $0<|\varepsilon|<<1$ and the $a_i$, $i=0,
1, 2, \ldots, k-1$, are real constants to be chosen in such a way
that the Poincar\'e map $\Pi(r_0; \varepsilon)$ associated to the
ordinary differential equation
\[ \frac{dr}{d\theta} \, = \, \frac{r^d \, R(r, \theta) \, + \,
\sum_{i=0}^{k-1} \varepsilon^{k-i} \, a_i \, r^{2i+d}}{r^{d-1} \,
F(r, \theta)} \, = \, \frac{r  R(r, \theta) \, + \, r \,
\sum_{i=0}^{k-1} \varepsilon^{k-i} \, a_i \, r^{2i}}{ F(r,
\theta)}  \] has $2k+1$ real zeroes; $k$ of them positive. We
recall that $F(0, \theta) \neq 0$ for all $\theta \in [0,2\pi)$
and, thus, the perturbative terms are an analytic perturbation in
a neighborhood of $r=0$ and $\varepsilon=0$ of equation
(\ref{eq5b}). The proof of the fact that this choice of $a_i$ can
be done is analogous to the one described in \cite{Andronov}, pp
254--259. More concretely, the exponent of the leading term of the
displacement function $d(r_0;0)$ of system (\ref{eq5b}) is $m$ and
the considered perturbation (\ref{eq5p}) produces that
$d(r_0;\varepsilon)$ has all the monomials of odd powers of $r_0$
up to order $m$. The coefficient of each monomial, for
$\varepsilon$ sufficiently small, is dominated by one of the
constants $a_i$. \par Undoing the change to polar coordinates,
system (\ref{eq5p}) gives rise to an analytic system in a
neighborhood of the origin which is a perturbation of system
(\ref{eqdeg}) of the form (\ref{eq2}) and with $k=(m-1)/2$ limit
cycles bifurcating from the origin. If system (\ref{eqdeg}) is
written as $\dot{x} \, = \, P(x,y)$ and $\dot{y} \, = \, Q(x,y)$,
then the change to cartesian coordinates from (\ref{eq5p}) reads
for: \begin{equation} \label{eqdegp1} \dot{x} \, = \, P(x,y) \, +
\, x \, K(x,y,\varepsilon), \qquad \dot{y} \, = \, Q(x,y)  \, + \,
y \, K(x,y,\varepsilon), \end{equation} where $\displaystyle
K(x,y,\varepsilon) \, = \, \sum_{i=0}^{k-1} \varepsilon^{k-i} \,
a_i \, (x^2+y^2)^{i+ \frac{d-1}{2}}$. We recall that $d$ is odd
and $d\geq 1$. In this way, we have that ${\rm
Cycl}(\mathcal{X}_\varepsilon,p_0) \geq k$. \par We provide now an
example of an analytic perturbation of system (\ref{eqdeg}) with
at least $(m+d)/2-1$ limit cycles bifurcating from the origin. We
take system (\ref{eqdegp1}) and we perturb it in order to produce
$\ell=(d-1)/2$ additional limit cycles bifurcating from the origin
when $\varepsilon \to 0$. Let us consider the smallest limit cycle
$\gamma$ surrounding the origin of system (\ref{eqdegp1}) and let
us assume that it is an attractor. We have that the $(m-1)/2$
limit cycles of system (\ref{eqdegp1}) which bifurcate from the
origin are hyperbolic, by choosing the parameters $a_i$
conveniently. Since $\gamma$ is the smallest limit cycle, we have
that the origin needs to be a repeller. Let us take a convenient
real value $b_{\ell-1}$ such that the system
\[ \begin{array}{lll} \displaystyle
\dot{x} & = &  \displaystyle P(x,y) \, + \, x \,
K(x,y,\varepsilon) \, + \,
\varepsilon^{k+1}\, x \, b_{\ell-1} \, (x^2+y^2)^{\ell-1}, \vspace{0.2cm} \\
 \displaystyle \dot{y} & = &  \displaystyle  Q(x,y) \, + \, y \, K(x,y,\varepsilon) \, + \,
\varepsilon^{k+1}\, y \, b_{\ell-1} \, (x^2+y^2)^{\ell-1}
\end{array}\] still has the previous limit cycle $\gamma$ as an
attractor (thus, $|b_{\ell-1}|$ needs to be small enough) and the
origin also becomes an attractor (this implies that
$b_{\ell-1}<0$). Therefore, there is a limit cycle bifurcating
from the origin, surrounding it and smaller than $\gamma$. If
$\gamma$ was a repeller and the origin an attractor in system
(\ref{eqdegp1}), we take $b_{\ell-1}>0$ in order to make the new
limit cycle to bifurcate. This bifurcated limit cycle is
hyperbolic by conveniently choosing the value $b_{\ell-1}$. The
previous system maintains the $(m-1)/2$ limit cycles of
(\ref{eqdegp1}) because they are all hyperbolic.  Thus, the
previous system has at least $(m-1)/2 +1$ limit cycles bifurcating
from the origin when $\varepsilon \to 0$.\par By induction, and a
relabelling of the parameters $a_i \, := \, b_{i+(d-1)/2}$, we
deduce that the following system
\begin{equation} \label{eqdegp2} \dot{x} \, = \, P(x,y) \, + \, x
\, \bar{K}(x,y,\varepsilon), \qquad \dot{y} \, = \, Q(x,y)  \, +
\, y \, \bar{K}(x,y,\varepsilon),
\end{equation} where $\displaystyle \bar{K}(x,y,\varepsilon) \, = \,
\sum_{i=0}^{L-1} \varepsilon^{L-i} \, b_i \, (x^2+y^2)^{i}$, $L \,
: =\, (m+d)/2-1$, has at least $(m+d)/2-1$ limit cycles
bifurcating from the origin. We recall that both $m$ and $d$ are
odd and $d\geq 1$, $m \geq 1$. In this way, we have that ${\rm
Cycl}(\mathcal{X}_\varepsilon,p_0) \geq (m+d)/2-1$. \newline

Finally, we will prove statement {\rm 2.} If we only consider
perturbations of system {\rm (\ref{eqdeg})} whose subdegree is
$\geq d$, that is vector fields of the type
$\mathcal{X}_\varepsilon^{[d]}$, then, the transformation of these
perturbative terms to polar coordinates gives rise to a
perturbation of the corresponding equation (\ref{eq5b}) which is
analytic in a neighborhood of $r=0$ and $\varepsilon=0$. Let us
assume that the circle $r=0$ is a limit cycle with multiplicity
$m$ and we consider the Poincar\'e map $\Pi(r_0)$ defined in the
previous section, which satisfies:
\[ \Pi(r_0) \, = \, r_0 \, + \, c_m r_0^m \, + \,
\mathcal{O}(r_0^{m+1}), \] with $c_m \neq 0$. We recall that
equation (\ref{eq5b}) has the discrete symmetry $(r, \theta)
\mapsto (-r, \theta+\pi)$  because it comes from a system in
cartesian coordinates (\ref{eqdeg}). This symmetry implies that
$r=0$ is a limit cycle which cannot be semistable and therefore
$m$ is odd. \par The key point of the proof is that any
perturbation of (\ref{eq5b}) analytic near $(r,
\varepsilon)=(0,0)$, with $\varepsilon \in \mathbb{R}^p$ small,
has a displacement function $d(r_0; \varepsilon) = \Pi(r_0;
\varepsilon)-r_0$ which is analytic near $(r_0,
\varepsilon)=(0,0)$ and when $\varepsilon=0$ coincides with the
displacement function $d(r_0; 0) \, = \, \Pi(r_0)-r_0$ of the
unperturbed equation (\ref{eq5b}). By using standard arguments
(counting zeroes with Weierstrass Preparation Theorem of $d(r_0;
\varepsilon)$ near $(r_0, \varepsilon)=(0,0)$), the cyclicity of
the circle $r=0$ under analytic perturbations of equation
(\ref{eq5b}) is $m$. We recall that the cyclicity and the
multiplicity of a limit cycle are equal, see \cite{Andronov}.
However, since the displacement function $d(r_0; 0) \, = \,
\Pi(r_0)-r_0$ of the unperturbed equation (\ref{eq5b}) is of odd
order $m$ at $r_0=0$, and taking into account the above mentioned
discrete symmetry, we have that only $(m-1)/2$ zeroes of $d(r_0;
\varepsilon)$ can appear for $r_0>0$ and $\| \varepsilon\|$ small
enough. This fact gives that at most $(m-1)/2$ limit cycles
bifurcate from the origin $p_0$ of system (\ref{eqdeg}) when only
this kind of perturbative terms are taken into account. Therefore,
we have proved that ${\rm Cycl}
(\mathcal{X}_\varepsilon^{[d]},p_0)=(m-1)/2$. Indeed, the example
given in (\ref{eqdegp1}) shows that this upper bound is sharp.
\bbox

\subsection*{Proof of Theorem {\rm \ref{thnil}}.}

The proof of this theorem goes analogously to the proof of Theorem
\ref{thdeg}, only with some technical differences. \newline

We first show that the transformation to generalized polar
coordinates (\ref{change}) $ x \, = \, r \, {\rm Cs} \, \theta, $
$ y \, = \, r^n \, {\rm Sn} \, \theta $ transforms system
(\ref{eqnil2}) to an equation over the cylinder of the form
(\ref{eq3}).
\begin{lemma} \label{lemnil1}
We assume that the origin of system {\rm (\ref{eqnil2})} is a
nilpotent monodromic singular point. Then the transformation to
generalized polar coordinates $\ (x,y) \, = \, (r \, {\rm Cs} \,
\theta, \, r^n \, {\rm Sn} \, \theta)$ brings system {\rm
(\ref{eqnil2})} to an ordinary differential equation {\rm
(\ref{eq3})} over a cylinder.
\end{lemma}
{\em Proof.} Taking into account that $\  {\rm Cs}^{2n} \theta \, +\,  n \, {\rm
Sn}^2 \theta\, =\, 1 $, we get that the Jacobian determinant of
the former change is
$$
J(r,\theta) = \frac{\partial (x,y)}{\partial (r,\theta)} = r^n .
$$
Since $x^{2n}+n y^{2} = r^{2n}$, we deduce that
\[
\dot{r} \, = \, \frac{x^{2n-1} \dot{x} + y \dot{y}}{r^{2n-1}} , \qquad
\dot{\theta} \, = \, \frac{x \dot{y} -n y \dot{x}}{r^{n+1}} .
\]
In particular, system (\ref{eqnil2}) adopts the form
\begin{equation}\label{eqnil3a}
\dot{r} \, = \,  \displaystyle \tilde{p}(\theta) \, r^{n+1} +
O(r^{n+2}) , \quad \dot{\theta}\,  = \,  \displaystyle r^{n-1} \ +
O(r^{n}) ,
\end{equation}
when $\beta > n-1 \ \mbox{or} \ \phi(x) \equiv 0$, and
\begin{equation}\label{eqnil3b}
\dot{r} \, = \, \displaystyle b \, {\rm Cs}^{n-1} \theta \, {\rm
Sn}^2\theta \, r^{n} + O(r^{n+1}) , \quad \dot{\theta} \, = \,
\displaystyle \left(1 + b \, {\rm Cs}^{n}\theta\, {\rm Sn}\,
\theta \right) r^{n-1} + O(r^{n}) ,
\end{equation}
when $\beta = n-1$. \par We also observe that in the latter case
($\beta = n-1$) we have the following decomposition \[ 1 + b \,
{\rm Cs}^{n}\theta\, {\rm Sn}\, \theta \, = \, \left({\rm
Cs}^{n}\theta\,  \, + \, \frac{b}{2}\, {\rm Sn}\, \theta \right)^2
\, + \, \frac{1}{4}\, (4n-b^2)\, {\rm Sn}^2\theta, \] where we
have used that $\  {\rm Cs}^{2n}\theta \, +\,  n \, {\rm Sn}^2
\theta\, =\, 1$. Since $4n-b^2>0$ in this case, due to Andreev's
conditions for monodromy, we deduce that $1 + b \, {\rm
Cs}^{n}\theta\, {\rm Sn}\, \theta>0$ for any $\theta \in
\mathbb{R}$. \par We denote by $\Xi(r, \theta)$ the function
defined by $\dot{r}$ and by $\Theta(r, \theta)$ the function
defined by $\dot{\theta}$ in both cases, and we have that,
\begin{equation} \label{eqnb}
\dot{r} \, = \, \Xi(r,\theta), \quad \dot{\theta} \, = \,
\Theta(r, \theta) \, = \, \Theta_{n-1}(\theta)\, r^{n-1} \, + \,
\mathcal{O}\left( r^n \right), \end{equation} where
\[ \Theta_{n-1}(\theta) \, = \, \left\{ \begin{array}{lll} \displaystyle 1 & \ \mbox{if} & \beta > n-1 \ \ \mbox{or} \ \
\phi(x) \equiv 0, \vspace{0.2cm} \\ \displaystyle 1 \, + \, b \,
{\rm Cs}^{n}\theta\, {\rm Sn}\, \theta  & \ \mbox{if} & \beta =
n-1.
\end{array} \right. \]
Hence, the equation of the orbits corresponding to system
(\ref{eqnil3a}) or (\ref{eqnil3b}) writes as
$$
\frac{d r}{d \theta} \, = \, \left\{ \begin{array}{lll}
{\displaystyle \frac{O(r^2)}{1+O(r)} } & \ \mbox{if} & \beta > n-1
\ \ \mbox{or} \ \ \phi(x) \equiv 0, \vspace{0.2cm} \\
{\displaystyle \frac{O(r)}{1+b \, {\rm Cs}^{n}\theta\, {\rm Sn}\,
\theta +O(r)} } & \ \mbox{if} & \beta = n-1.
\end{array} \right.
$$
We observe that $\Theta_{n-1}(\theta) \neq 0$ for any $\theta \in
[0,T_n)$. In short, we have proved that in a neighborhood of any
monodromic singular point of the form (\ref{eqnil}), we can
perform a transformation, which is the composition of the changes
(\ref{change1}) and (\ref{change2}) and the transformation to
generalized polar coordinates, which brings the system to an
equation over a cylinder $C$:
\begin{equation} \label{eqnp} \frac{d r}{d \theta} \, = \, \mathcal{F}(r, \theta) , \end{equation}
where $\mathcal{F}(r,\theta) $ is $T_n$--periodic in $\theta$ and
$\mathcal{F}(0,\theta) \equiv 0$. \bbox \newline

The center problem for the origin of system (\ref{eqnil}) is
equivalent to determine when the circle $r=0$ is contained in a
period annulus for equation (\ref{eqnp}).

\begin{lemma} \label{lemnil1b}
We assume that the origin of system {\rm (\ref{eqnil})} is a
monodromic singular point with Andreev number $n$. The circle
$r=0$ is a limit cycle of equation {\rm (\ref{eqnp})} if, and only
if, the origin of system {\rm (\ref{eqnil})} is a focus.
\end{lemma}
{\em Proof.} The proof is analogous to the proof of Lemma
\ref{lemdeg1b}. The transformation from system (\ref{eqnil}) to
equation (\ref{eqnp}) gives a one-to-one correspondence between
each point in a punctured neighborhood of the origin in the plane
$(x,y)$ and each point on a cylinder $\{ (r,\theta) \, : \,
0<r<\delta, \theta \in \mathcal{S}^1\}$ for $\delta>0$
sufficiently small. \bbox \newline

The following proposition establishes a symmetry for equation
(\ref{eqnp}) which is inherited by the symmetries of the
generalized trigonometric functions.
\begin{proposition} \label{propnil}
Let us consider a planar $\, \mathcal{C}^1$ differential system
$\dot{x} \, = \, P(x,y)$, $\dot{y} \, = \, Q(x,y)$, we take any
positive integer $n$ and perform the change to generalized polar
coordinates $x \, = \, r \, {\rm Cs} \, \theta$, $y \, = \, r^n \,
{\rm Sn} \, \theta$. The resulting system $\dot{r} \, = \, \Xi(r,
\theta)$, $\dot{\theta} \, = \, \Theta(r, \theta)$ is invariant
under the change of variables $ (r, \theta) \, \mapsto \,
\left(-r, (-1)^{n+1} \, \left[\theta + T_n/2 \right]\right). $
\end{proposition}
{\it Proof.} We observe that, due to Proposition \ref{propgen},
the following composition is the identity ($X=x$, $Y=y$):
\[ (x, y) \, \mapsto \, (r, \theta) \, \mapsto \, (R, \varphi) \, \mapsto \, (X,Y), \]
where  $R=-r$, $\varphi =  (-1)^{n+1} [\theta + T_n/2]$ and $X \,
= \, R \, {\rm Cs} \, \varphi$, $Y \, = \, R^n \, {\rm Sn} \,
\varphi$. Since \[ \dot{r} \, = \, \frac{x^{2n-1} P(x,y) + y
Q(x,y)}{r^{2n-1}} , \qquad \dot{\theta} \, = \, \frac{x Q(x,y) -n
y P(x,y)}{r^{n+1}} ,
\]
the proposition follows.  \bbox \newline

The previous symmetry of equation (\ref{eqnp}) imposes a condition
on the circle $r=0$ to be a limit cycle.

\begin{lemma} \label{lemdnil}
We assume that the origin of system {\rm (\ref{eqnil})} is a focus
with Andreev number $n$. We consider the corresponding equation
{\rm (\ref{eqnp})} and its Poincar\'e map $\Pi(r_0) \, = \, r_0 \,
+ \, c_m r_0^m \, + \, \mathcal{O}(r_0^{m+1}),$ with $c_m \neq 0$.
\begin{itemize}
\item[{\rm (i)}] If $n$ is odd, then $r=0$ cannot be a semistable limit cycle of equation {\rm (\ref{eqnp})}, that is, $m$ is odd.
\item[{\rm (ii)}] If $n$ is even, then $r=0$ is a semistable limit cycle of equation {\rm (\ref{eqnp})}, that is, $m$ is even.
\end{itemize}  \end{lemma}
{\em Proof.} The change to generalized polar coordinates ensures
that the origin of system (\ref{eqnil}) is a focus if, and only
if, the circle $r=0$ is a limit cycle of equation {\rm
(\ref{eqnp})}. We assume that the circle $r=0$ is a limit cycle
with multiplicity $m$ and we consider its associated Poincar\'e
map $\Pi(r_0)$. We recall that equation (\ref{eqnp}) has the
discrete symmetry $(r, \theta) \mapsto (-r, (-1)^{n+1} \,
\left[\theta + T_n/2 \right])$  because it comes from a system in
cartesian coordinates (\ref{eqnil}), see Proposition
\ref{propnil}. This symmetry implies that $r=0$ is a semistable
limit cycle if, and only if, $n$ is even. \bbox

\begin{corollary} \label{nilp20}
If equation {\rm (\ref{eqnp})} has a periodic orbit different from
$r = 0$, then it has two periodic orbits (one in the upper half
cylinder and one in the lower half cylinder).
\end{corollary}
{\it Proof.} The properties of the periodic orbits of (\ref{eqnp})
stated in this corollary are straightforward consequences of the
discrete symmetry given in Proposition \ref{propnil}. \bbox
\newline

Let $V_0(x,y)$ be an inverse integrating factor defined in a neighborhood of the
origin for system (\ref{eqnil}). Since the changes of variables
(\ref{change1}) and (\ref{change2}) have constant
non--vanishing Jacobian, it follows that the transformed system
(\ref{eqnil2}) has the following inverse integrating
factor defined in a neighborhood of the
origin:
\[
V_0^*(x,y) \, =\,  V_0(\xi^{-1} x ,  -\xi^{-1} y + F(\xi^{-1} x)) .
\]
\par

Next lemma is the analogous to Lemma \ref{lemdeg2} and states
that, under our hypothesis, we have an inverse integrating factor
$V(r, \theta)$ for equation (\ref{eqnp}) which is analytic in $r$
in a neighborhood of $r=0$ and such that $V(0, \theta) \equiv 0$.

\begin{lemma} \label{lemnil2} We assume that the origin of system {\rm (\ref{eqnil2})} is a
nilpotent monodromic singularity.
\begin{itemize}
\item[{\rm (i)}] If system {\rm (\ref{eqnil2})} has an inverse integrating
factor $V_0^*(x,y)$ defined in a neighborhood of the origin, then the function defined by
\[ V(r, \theta) \, := \, \frac{V_0^*(r \, {\rm Cs} \, \theta, r^n {\rm
Sn}\, \theta)}{r^{n} \, \Theta(r, \theta)} , \] where $\Theta(r,
\theta)$ is the function defined in {\rm (\ref{eqnb})}, is an
inverse integrating factor for equation {\rm (\ref{eqnp})} in $r
\neq 0$.
\item[{\rm (ii)}] Let $V(r,\theta)$ be an inverse integrating factor of equation {\rm (\ref{eqnp})}. We assume that $V(r,  \theta )$ has a Laurent expansion in a neighborhood of $r=0$ of the form $ V(r, \theta) \, = \, v_{m}(\theta) \, r^m \, + \, \mathcal{O}(r^{m+1}), $ with $v_m(\theta) \not\equiv 0$ and $m \in \mathbb{Z}$. Then, if $m \leq 0$ the origin of system
{\rm
(\ref{eqnil2})} is a center.
\end{itemize}
\end{lemma}
{\em Proof.} {\rm (i)} Taking into account the Jacobian $r^n$ of
the change to generalized polar coordinates $\ (x,y) \, = \, (r \,
{\rm Cs} \, \theta, \, r^n \, {\rm Sn} \, \theta)$, we see that
the differential equation (\ref{eqnp}) has the inverse integrating
factor $V(r,\theta)$ described in the statement which is a
$T_n$--periodic function of $\theta$. We observe that
$V(r,\theta)$ may not be well-defined on $r=0$. \par {\rm (ii)}
The assumption $m \leq 0$ establishes that we are under the
hypothesis  of Lemma \ref{lemv0} and we conclude that the cycle
$r=0$ is surrounded by periodic orbits, which give rise to a
neighborhood of the origin of system (\ref{eqnil2}) filled with
periodic orbits. \bbox
\newline

The previous Lemmas \ref{lemdnil} and \ref{lemnil2}, together with
Corollary \ref{corvpi}, ensure that, under the hypothesis of
Theorem \ref{thnil} and if the origin of system (\ref{eqnil}) is a
focus, then $m \geq 1$ and $m+n$ needs to be even. Thus, we have
proved statement (i) in Theorem \ref{thnil}.
\newline

To end with, we relate the multiplicity of the limit cycle $r=0$
of equation (\ref{eqnp}) and the cyclicity of the origin of system
(\ref{eqnil}).

\begin{lemma}
We assume that the origin of system {\rm (\ref{eqnil})} is a
nilpotent focus with Andreev number $n$. We consider the
corresponding equation {\rm (\ref{eqnp})} for which we assume the
circle $r=0$ to be a limit cycle with multiplicity $m$.
\begin{itemize}
\item[{\rm 1.}] The cyclicity ${\rm
Cycl}(\mathcal{X}_\varepsilon,p_0)$ of the origin of system {\rm
(\ref{eqnil})} satisfies ${\rm Cycl}(\mathcal{X}_\varepsilon,p_0)
\geq (m+n)/2-1$.
\item[{\rm 2.}] If only analytic perturbations of system {\rm (\ref{eqnil2})} with $(1,n)$--quasihomogeneous weighted subdegrees
$(w_x,w_y)$ with $w_x \geq n$ and $w_y \geq 2n-1$ are taken into
account, then the maximum number of limit cycles which bifurcate
from the origin of system {\rm (\ref{eqnil2})} {\rm (}and, thus,
of system {\rm (\ref{eqnil})}{\rm )} is $\lfloor (m-1)/2 \rfloor$,
that is, ${\rm Cycl}(\mathcal{X}_\varepsilon^{[n,2n-1]},p_0) \, =
\, \lfloor (m-1)/2 \rfloor$.
\end{itemize}
\end{lemma}
{\em Proof.} The fact that $m$ and $n$ have the same parity is a
consequence of Lemma \ref{lemdnil}. If $n$ is odd, the symmetry
given in Proposition \ref{propnil} implies that at most $(m-1)/2$
limit cycles can bifurcate from the limit cycle $r=0$ of equation
(\ref{eqnp}) in the region $r>0$. If $n$ is even, since equation
(\ref{eqnp}) comes from cartesian coordinates, we also need to
take into account that $r=0$ is always a solution. Therefore, by
using the symmetry again, at most $(m-2)/2$ limit cycles can
bifurcate from $r=0$ in the region $r>0$.
\par We provide an example of a perturbation of equation
(\ref{eqnp}), with $m$ limit cycles bifurcating from $r=0$
(counting multiplicities), whose transformation to cartesian
coordinates gives a perturbation of system (\ref{eqnil}) with
exactly $\lfloor (m-1)/2 \rfloor$ limit cycles bifurcating from
the origin. This example proves that ${\rm
Cycl}(\mathcal{X}_\varepsilon,p_0) \geq k$. We take equation
(\ref{eqnp}) and we assume that $r=0$ is a limit cycle of
multiplicity $m$, which is an integer with the same parity as $n$
and such that $m \geq 1$. We consider the associated system
(\ref{eqnb}) from which (\ref{eqnp}) comes from, we define
$k=\lfloor (m-1)/2 \rfloor$ and we perturb system (\ref{eqnb}) in
the following way:
\begin{equation}
\begin{array}{lll}
\displaystyle \dot{r} \, = \, \Xi(r, \theta) \, + \,
\sum_{i=0}^{k-1} \varepsilon^{k-i} \, a_i \, r^{n+2i} \left( {\rm
Cs} \, \theta\right)^{n-1+2i}, & \ \dot{\theta} \, = \, \Theta(r,
\theta), & \ \mbox{if $n$ is odd}, \vspace{0.2cm} \\
\displaystyle \dot{r} \, = \, \Xi(r, \theta) \, + \,
\sum_{i=0}^{k-1} \varepsilon^{k-i} \, a_i \, r^{n+1+2i} \left(
{\rm Cs} \, \theta\right)^{n+2i}  , & \ \dot{\theta} \, = \,
\Theta(r, \theta), & \ \mbox{if $n$ is even}, \end{array}
\label{eqpp}
\end{equation} with the convention that if $k=0$ no perturbation
terms are taken. The real value $\varepsilon$ is the perturbation
parameter $0<|\varepsilon|<<1$ and the $a_i$, $i=0, 1, 2, \ldots,
k-1$, are real constants to be chosen so that the Poincar\'e map
$\Pi(r_0; \varepsilon)$ associated to the ordinary differential
equation \[ \frac{dr}{d\theta} \, = \, \left\{ \begin{array}{ll}
\displaystyle \, \frac{\Xi(r, \theta) \, + \, \sum_{i=0}^{k-1}
\varepsilon^{k-i} \, a_i \, r^{n+2i} \left( {\rm
Cs} \, \theta\right)^{n-1+2i}}{\Theta(r, \theta)}, & \ \mbox{if $n$ is odd}, \vspace{0.2cm} \\
\displaystyle \, \frac{\Xi(r, \theta) \, + \, \sum_{i=0}^{k-1}
\varepsilon^{k-i} \, a_i \,  r^{n+1+2i} \left( {\rm Cs} \,
\theta\right)^{n+2i} }{\Theta(r, \theta)}, & \ \mbox{if $n$ is
even}, \end{array} \right. \] has $2k+1$ real zeroes; $k$ of them
positive. The proof of this fact is analogous to the proof given
in \cite{Andronov}, pp. 254--259.
\par Undoing the change to generalized polar coordinates system
(\ref{eqpp}) gives rise to an analytic system in a neighborhood of
the origin which is a perturbed system from (\ref{eqnil2}) of the
form (\ref{eq2}) and with $k=\lfloor (m-1)/2 \rfloor$ limit cycles
bifurcating from the origin. All these limit cycles are hyperbolic
by taking convenient values of the parameters $a_i$. If system
(\ref{eqnil2}) is written as $\dot{x} \, = \, P(x,y)$ and $\dot{y}
\, = \, Q(x,y)$, then the change to cartesian coordinates from
(\ref{eqpp}) reads for:
\begin{equation} \label{eqnilp1} \dot{x} \, = \, P(x,y) \, + \, x \, K(x,y,\varepsilon),
\qquad \dot{y} \, = \, Q(x,y) + n y K(x,y,  \varepsilon),
\end{equation} where
\[ K(x,y , \varepsilon) \, = \, \left\{ \begin{array}{ll} \displaystyle \sum_{i=0}^{k-1}
\varepsilon^{k-i} \, a_i \, x^{n-1+2i} & \ \mbox{if $n$ is odd},
\vspace{0.2cm} \\
\displaystyle \sum_{i=0}^{k-1} \varepsilon^{k-i} a_i \, x^{n+2i} &
\ \mbox{if $n$ is even}. \end{array} \right. \] Then, undoing the
changes (\ref{change1}) and (\ref{change2}), we obtain a
perturbation of system (\ref{eqnil}) with $k=\lfloor (m-1)/2
\rfloor$ limit cycles bifurcating from the origin. We observe that
these perturbative terms satisfy the $(1,n)$-quasihomogeneous
subdegree conditions established in the statement 2 of the lemma.
\par In order to generate the $\ell \, := \, \lfloor n/2 \rfloor$
limit cycles that we lack to prove the first statement of this
lemma, we will consider perturbations with lower subdegree. We
observe that the stability of the origin in system (\ref{eqnilp1})
is given by the sign of $a_0$. Thus, if we choose a convenient
value $b_{\ell-1}$ such that the following system
\[ \begin{array}{lll} \displaystyle \dot{x} & = & \displaystyle
P(x,y) \, + \, x \, K(x,y,\varepsilon) \, + \, \varepsilon^{k+1}
\, x \, b_{\ell-1} \, x^{2(\ell-1)}, \vspace{0.2cm} \\
 \displaystyle \dot{y} & = & Q(x,y) + n y K(x,y,\varepsilon) \, + \, \varepsilon^{k+1}
\, y \, b_{\ell-1} \, x^{2(\ell-1)} \end{array} \] satisfies that
the stability of the smallest limit cycle in (\ref{eqnilp1}) does
not change (thus $|b_{\ell-1}|$ needs to be small enough) and with
the origin of the contrary stability (namely, $b_{\ell-1}$ of
contrary sign to $a_0$), then there is a new limit cycle
bifurcating from the origin. This limit cycle is hyperbolic by
choosing $b_{\ell-1}$ conveniently. By induction, and relabelling
$a_i \, := b_{i+\ell}$, we have that the system
\begin{equation} \label{eqnilp2} \dot{x} \, = \, P(x,y) \, + \, x
\, \bar{K}(x,y,\varepsilon), \qquad \dot{y} \, = \, Q(x,y) + n y
\bar{K}(x,y, \varepsilon),
\end{equation} where \[ \bar{K}(x,y , \varepsilon) \, = \, \displaystyle \sum_{i=0}^{L-1}
\varepsilon^{L-i} \, b_i \, x^{2i} \] and $L=(m+n)/2-1$, has at
least $(m+n)/2-1$ limit cycles bifurcating from the origin. By
undoing the changes (\ref{change1}) and (\ref{change2}), we have
thus proved that the cyclicity of the origin of system
(\ref{eqnil}) is at least $(m+n)/2-1$. \newline

We are going to prove statement 2 of this lemma. By the results
established on Section \ref{sect2}, we can control the maximum
number of limit cycles which bifurcate from $r=0$ in equation
(\ref{eqnp}) under analytic perturbations of this equation. The
hypothesis that we only take analytic perturbations of system {\rm
(\ref{eqnil2})} with $(1,n)$--quasihomogeneous weighted subdegrees
$(w_x,w_y)$ with $w_x \geq n$ and $w_y \geq 2n-1$ is equivalent to
say that we take a perturbation of equation (\ref{eqnp}) which is
analytic in a neighborhood of both $r=0$ and $\varepsilon=0$.
System (\ref{eqnilp1}) provides an example where the upper bound
of $k\, = \, \lfloor (m-1)/2 \rfloor$ limit cycles bifurcating
from the origin $p_0$ of system (\ref{eqnil}), when only these
perturbations are considered, is attained. In this way we have
that ${\rm Cycl}(\mathcal{X}_\varepsilon^{[n,2n-1]},p_0) \, = \,
\lfloor (m-1)/2 \rfloor$. \bbox

\subsubsection*{Proof of Corollary \ref{cornil}.}

Let us consider system (\ref{eqnil}) with a nilpotent focus at the
origin and the corresponding Andreev number $n$. Let $V_0(x,y)$ be
an inverse integrating factor of this system which is analytic in
a neighborhood of the origin. Therefore, we have an inverse
integrating factor $V_0^{*}(x,y)$ of the corresponding system
(\ref{eqnil2}) which is analytic in a neighborhood of the origin.
We consider the Taylor development around $r=0$ of the following
function:
\[ V_0^{*}(r \, {\rm Cs} \, \theta, \, r^n \, {\rm Sn}\, \theta)
\, = \, v_{M}^{*}(\theta) \, r^M \, + \, \mathcal{O}(r^{M+1}), \]
where $v_{M}^{*}(\theta) \not\equiv 0$. \par Let us consider the
transformation of system (\ref{eqnil2}) to an equation over a
cylinder by generalized polar coordinates, see Lemma
\ref{lemnil1}, and the corresponding inverse integrating factor
$V(r,\theta)$ which is smooth and non-flat in $r$ in a
neighborhood of $r=0$. We take the Taylor development around $r=0$
of the following functions:
\[ V(r,\theta) \, = \, v_m(\theta) \, r^m \, + \,
\mathcal{O}(r^{m+1}), \quad \Theta(r,\theta) \, = \,
\Theta_{n-1}(\theta) r^{n-1} \, + \, \mathcal{O}(r^n), \] where
$\Theta(r,\theta)$ is the function defined in (\ref{eqnb}). We
recall, see again the proof of Lemma \ref{lemnil1}, that
$\Theta_{n-1}(\theta) \neq 0$ for any $\theta \in [0,T_n)$.
Indeed, by Lemma \ref{lemv1} we have that $v_{m}(\theta) \neq 0$
for any value of $\theta \in [0,T_n)$. \par By statement (i) in
Lemma \ref{lemnil2} we deduce that $v_{M}^{*}(\theta) \, = \,
v_m(\theta) \, \Theta_{n-1}(\theta)$. Therefore, we have that
$v_{M}^{*}(\theta) \neq 0$ for any value of $\theta \in [0,T_n)$.
On the other hand, since $V_0^{*}(x,y)$ is analytic in a
neighborhood of the origin, we deduce that $v_{M}^{*}(\theta)$ is
a $(1,n)$--quasihomogeneous trigonometric polynomial of weighted
degree $M$. Using the symmetries of the generalized trigonometric
functions described in statement {\rm (e)} of Proposition
\ref{propgen}, we see that \[ v_M^{*} \left( \frac{T_n}{2} \right)
\, = \, (-1)^M \, v_M^{*}(0).
\] Therefore, we conclude
that $M$ needs to be an even number.
\par We consider the value of $m$ defined in Theorem \ref{thnil} and Remark
\ref{remnil}, and we have that $M=2n-1+m$. Since $M$ is even, we
deduce that $m$ is odd. Moreover, since the origin of system
(\ref{eqnil}) is assumed to be a focus, we have that $m$ is an
integer number with $m\geq 1$ and that $m$ and $n$ need to have
the same parity, see statement {\rm (ii)} of Theorem \ref{thnil}.
Thus, $n$ is odd. \bbox

\subsubsection*{Proof of Corollary \ref{corexist}.}

As we have already stated, see Lemmas \ref{lemdeg1} and
\ref{lemnil1}, in a neighborhood of these singular points, we can
transform the system to an equation over a cylinder by means of
(generalized) polar coordinates. Since the origin is a focus, we
have that the periodic orbit $r=0$ is a limit cycle for the
equation on the cylinder. Lemma \ref{lemv2} ensures the existence
of a smooth, and non--flat in $r$ in a neighborhood of $r=0$,
inverse integrating factor $V(r,\theta)$ for the equation over the
cylinder. Indeed, we have that its vanishing multiplicity $m$ at
the origin is at least $1$, that is, there exists a smooth
function $f(r,\theta)$ defined on the cylinder such that
$V(r,\theta) \, = \, r^m\, f(r,\theta)$, with $m \geq 1$. The
inverse integrating factor $V(r, \theta)$ gives rise to an inverse
integrating factor $V_0(x,y)$ in cartesian coordinates. \par In
the case of polar coordinates, we have that $r=\sqrt{x^2+y^2}$ and
$\theta=\arctan(y/x)$. We observe that the functions
$\frac{\partial r}{\partial x}$, $\frac{\partial r}{\partial y}$,
$r \frac{\partial \theta}{\partial x}$ and $r \frac{\partial
\theta}{\partial y}$ are bounded functions in a neighborhood of
the origin. By Remark \ref{remdeg}, and if we consider system
(\ref{eqdeg}) with $d \geq 1$, we have that $V_0(x,y) \, = \,
r^{m+d} \, \tilde{f}(r,\theta)$ with $r$ and $\theta$ expressed in
cartesian coordinates and $\tilde{f}(r,\theta)$ a bounded
function, with bounded derivatives, in a neighborhood of the
origin.
\par  In the case of generalized polar coordinates, we have that $r \, = \,
\sqrt[2n]{x^{2n}+ny^2}$ and we observe that the functions
$\frac{\partial r}{\partial x}$ and $r^{n-1}\frac{\partial
r}{\partial y}$ are bounded in a neighborhood of the origin
because
\[ \frac{\partial r}{\partial x} \, = \,
\frac{x^{2n-1}}{(x^{2n}+ny^2)^{\frac{2n-1}{2n}}} \quad \mbox{and}
\quad \frac{\partial r}{\partial y} \, = \,
\frac{y}{(x^{2n}+ny^2)^{\frac{2n-1}{2n}}}. \] If, in these
expressions, we consider again the change of coordinates
(\ref{change}) we have that $\frac{\partial r}{\partial x} \, = \,
{\rm Cs}^{2n-1} \theta$ and $r^{n-1}\frac{\partial r}{\partial y}
\, = \, {\rm Sn}\, \theta$. $ $From the change of coordinates
(\ref{change}) and using the definition of the generalized
trigonometric functions, see (\ref{Cauchy}), it can be shown that
\[ \frac{\partial \theta}{\partial x} \, = \,
\frac{-\, n\, y}{(x^{2n}+ny^2)^{\frac{n+1}{2n}}} \quad \mbox{and}
\quad \frac{\partial \theta}{\partial y} \, = \,
\frac{x}{(x^{2n}+ny^2)^{\frac{n+1}{2n}}}. \] Thus, we have that
the functions $r \frac{\partial \theta}{\partial x}$ and $r^n
\frac{\partial \theta}{\partial y}$ are bounded in a neighborhood
of the origin by an analogous argument as before. By Remark
\ref{remnil}, and if we consider system (\ref{eqnil}) with $n >
1$, we have that $V_0(x,y) \, = \, r^{m+2n-1} \,
\tilde{f}(r,\theta)$ with $r$ and $\theta$ expressed in cartesian
coordinates and $\tilde{f}(r,\theta)$ a bounded function, with
bounded derivatives, in a neighborhood of the origin.
\par We have, in both cases, that $V_0(x,y) \, = \, r^{a} \,
\tilde{f}(r,\theta)$, with $r$ and $\theta$ expressed in cartesian
coordinates and $\tilde{f}(r,\theta)$ a bounded function, with
bounded derivatives, in a neighborhood of the origin. Since $m
\geq 1$, we have that the exponent $a>1$ in the case of polar
coordinates and $a>n$ in the case of generalized polar
coordinates. Thus, the limit of this product when $(x,y)$ tends to
the origin exists and it is equal to zero. Therefore, $V_0(x,y)$
is continuous in a neighborhood of the origin and $V_0(0,0)=0$.
\par By the chain rule, we have that
\[ \frac{\partial V_0}{\partial x} \, = \, \frac{\partial
V_0}{\partial r} \, \frac{\partial r}{\partial x} \, + \,
\frac{\partial V_0}{\partial \theta} \, \frac{\partial
\theta}{\partial x} \, = \, \left( a r^{a-1} \tilde{f}(r,\theta)
\, + \, r^a \, \frac{\partial \tilde{f}}{\partial r} \right) \,
\frac{\partial r}{\partial x} \, + \, r^a \, \frac{\partial
\tilde{f}}{\partial \theta} \, \frac{\partial \theta}{\partial x}.
\] This expression can also be written as the product of a
function tending to zero (because $a>1$ or $a>n$, respectively, in
each case) and a bounded function. Thus, the function
$\frac{\partial V_0}{\partial x}$ is continuous in a neighborhood
of the origin and it is zero on the point $(0,0)$. An analogous
argument holds for $\frac{\partial V_0}{\partial y}$. We have that
the function $V_0(x,y)$ and its first derivatives are continuous
and vanish at the origin. Therefore, $V_0(x,y)$ is at least of
class $\mathcal{C}^1$ in a neighborhood of the origin. \bbox
\newline

We remark that in Example \ref{ex4} we have given an inverse
integrating factor $V_0(x,y)$ for system (\ref{exen}) which might
be only of class $\mathcal{C}^1$ in a neighborhood of the origin,
depending on the values of $m$ and $n$.
\newline

\vspace{0.5cm}

{\bf Addresses and e-mails:} \\
$^{\ (1)}$ Departament de Matem\`atica. Universitat de Lleida.
\\ Avda. Jaume II, 69. 25001 Lleida, SPAIN.
\\ {\rm E--mails:} {\tt garcia@matematica.udl.cat}, {\tt
mtgrau@matematica.udl.cat} \vspace{0.2cm}
\\
$^{\ (2)}$ Laboratoire de Math\'ematiques et Physique Th\'eorique.
C.N.R.S. UMR 6083. \\ Facult\'e des Sciences et Techniques.
Universit\'e de Tours. \\ Parc de Grandmont 37200 Tours, FRANCE.
\\ E-mail: {\tt Hector.Giacomini@lmpt.univ-tours.fr}

\end{document}